\documentclass[a4paper,14pt]{article}

\usepackage[margin=1.5in]{geometry}

\usepackage[english]{babel}
\usepackage[T1]{fontenc}
\usepackage[utf8]{inputenc}
\usepackage{lmodern}
\usepackage{microtype}

\usepackage{bbm}
\usepackage{amsmath}
\usepackage{amssymb}
\usepackage{enumerate}
\usepackage{amsthm}

\usepackage{xcolor}

\usepackage{listings}

\usepackage{graphicx}
\usepackage{epsfig}

\usepackage{url}
\usepackage[backend=bibtex,citestyle=numeric]{biblatex}
\usepackage{hyperref}

\usepackage{comment} 

\title{The k-PDTM : a coreset for robust geometric inference\footnote{This work was partially supported by the ANR project TopData and GUDHI}}

\author{
  Br\'echeteau, Claire\\
  \texttt{claire.brecheteau@inria.fr}\\
  Universit\'e Paris-Saclay -- LMO \& Inria
  \and
  Levrard, Cl\'ement\\
  \texttt{levrard@math.univ-paris-diderot.fr}\\
  Universit\'e Paris Diderot -- LPMA
}

\newtheoremstyle{break}  
  {\topsep}   
  {\topsep}   
  {\itshape}  
  {0pt}       
  {\bfseries} 
  {.}         
  {\newline}  
  {}          

\newtheoremstyle{break2}  
  {\topsep}   
  {\topsep}   
  {\itshape}  
  {0pt}       
  {\bfseries} 
  {}         
  {\newline}  
  {}          

\theoremstyle{break}
\newtheorem{thmm}{Theorem}
\newtheorem{dff}[thmm]{Definition}
\newtheorem{propp}[thmm]{Proposition}
\newtheorem{corr}[thmm]{Corollary}
\newtheorem{lmm}[thmm]{Lemma}

\theoremstyle{break2}
\newtheorem{Preuve}{Proof}

\theoremstyle{plain}
\newtheorem{Remarque}{Remark}
\newtheorem{exx}[thmm]{Example}

\newenvironment{pvv}{\begin{Preuve}\rm}%
  {$\blacksquare$\end{Preuve}}
\newenvironment{rqq}{\begin{Remarque}\rm}%
  {\end{Remarque}}

\newcommand{\thm}{\begin{thmm}}
\newcommand{\ethm}{\end{thmm}}
\newcommand{\lm}{\begin{lmm}}
\newcommand{\elm}{\end{lmm}}
\newcommand{\ex}{\begin{exx}}
\newcommand{\eex}{\end{exx}}
\newcommand{\df}{\begin{dff}}
\newcommand{\edf}{\end{dff}}
\newcommand{\prop}{\begin{propp}}
\newcommand{\eprop}{\end{propp}}
\newcommand{\pv}{\begin{pvv}}
\newcommand{\epv}{\end{pvv}}
\newcommand{\rem}{\begin{rqq}}
\newcommand{\erem}{\end{rqq}}
\newcommand{\cor}{\begin{corr}}
\newcommand{\ecor}{\end{corr}}

\usepackage{bbm}
\usepackage{color, colortbl}

\lstdefinestyle{cust}{
language=python,
commentstyle=\ttfamily,
basicstyle=,
escapeinside={\%*}{*)},
frame=single,
keepspaces=true,
keywordstyle=\bfseries,
morekeywords={*,Input,Output},
}
\usepackage{xcolor}
\newcommand{\R}{\mathbbm R}
\newcommand{\N}{\mathbbm N}
\newcommand{\E}{\mathbbm E}
\newcommand{\X}{\mathbbm X}

\newcommand{\1}{\mathbbm 1}
\newcommand{\dd}{\rm d\it}
\newcommand{\NN}{\mathcal N}
\newcommand{\CC}{\mathcal C}
\newcommand{\MM}{\mathcal M}
\newcommand{\A}{\mathcal A}

\newcommand{\PPP}{\mathcal P}
\newcommand{\B}{\rm B\it}
\newcommand{\HP}{\rm H\it}
\newcommand{\Sph}{\rm S\it}
\newcommand{\supp}{\rm Supp \it}
\newcommand{\conv}{\rm Conv \it}
\newcommand{\opt}{\mathcal{O}pt}
\newcommand{\OPT}{\rm OPT\it}
\newcommand{\Nn}{\rm NN\it}
\newcommand{\NNN}{\mathcal{NN}}

\newcommand{\dR}{\rm d\it_{\overline{\R}_d}}

\DeclareMathOperator*{\argmin}{arg\,\min}
\definecolor{Gray}{gray}{0.9}

\bibliography{biblio}

\begin{document}

\maketitle

\begin{abstract}
Analyzing the sub-level sets of the distance to a compact sub-manifold of $\R_d$ is a common method in TDA to understand its topology. 
The distance to measure (DTM) was introduced by Chazal, Cohen-Steiner and M\'erigot in \cite{Merigot1} to face the non-robustness of the distance to a compact set to noise and outliers. This function makes possible the inference of the topology of a compact subset of $\R_d$ from a noisy cloud of $n$ points lying nearby in the Wasserstein sense. In practice, these sub-level sets may be computed using approximations of the DTM such as the $q$-witnessed distance \cite{Merigot2} or other power distance \cite{Buchet16}. These approaches lead eventually to compute the homology of unions of $n$ growing balls, that might become intractable whenever $n$ is large.

       To simultaneously face the two problems of large number of points and noise, we introduce the $k$-power distance to measure ($k$-PDTM). This new approximation of the distance to measure may be thought of as a $k$-coreset based approximation of  the DTM. Its sublevel sets consist in union of $k$-balls, $k <<n$, and this distance is also proved robust to noise. We assess the quality of this approximation  for $k$ possibly dramatically smaller than $n$, for instance $k=n^{\frac{1}{3}}$ is proved to be optimal for $2$-dimensional shapes. We also provide an algorithm to compute this $k$-PDTM.
\end{abstract}


\section{Introduction}
\subsection{Background on robust geometric inference}
       Let $M \subset \mathbb{R}^d$ be a compact set included in the closed Euclidean ball $\overline{B}(0,K)$, for $K>0$, whose topology is to be inferred. A common approach is to sample $\X_n=\{X_1,X_2,\ldots,X_n\} $ on $M$, and approximate the distance to $M$ via the distance to the sample points. As emphasized in \cite{Merigot1}, such an approach suffers from non-robustness to outliers. To face this issue, \cite{Merigot1} introduces the \textit{distance to measure} as a robust surrogate of the distance to $M$, when  $\X_n=\{X_1,X_2,\ldots,X_n\} $ is considered as a $n$-sample, that is $n$ independent realizations of a distribution measure $P$ whose support $Supp(P)$ is $M$, possibly corrupted by noise. Namely, for a Borel probability measure $P$ on $\R_d$, a mass parameter $h\in\left[0,1\right]$ and $x\in\R_d$, the distance of $x$ to the measure $P$, $\dd_{P,h}(x)$ is defined by 
\df[DTM]
\label{def:DTM}
\[\dd^2_{P,h}(x)=\frac{1}{h}\int_{l=0}^h\delta^2_{P,l}(x)\,\dd l, \qquad \mbox{with} \qquad \delta_{P,l}(x)=\inf\{r>0\mid P(\overline{\B}(x,r))>l\},\]
\edf
where $\overline{\B}(x,r)$ denotes the closed Euclidean ball with radius $r$. When $P$ is uniform enough on a compact set with positive reach $\rho$, this distance is proved to approximate well the distance to $M$ (\cite[Proposition 4.9]{Merigot1}) and is robust to noise (\cite[Theorem 3.5]{Merigot1}). The distance to measure is usually inferred from $\X_n$ via its empirical counterpart, also called \textit{empirical DTM}, replacing $P$ by the empirical distribution $P_n=\frac{1}{n}\sum_{i=1}^n\delta_{X_i}$, where $\delta_{x}$ is the Dirac mass on ${x}$.  

M\'erigot et al noted in \cite{Merigot2} that the sublevel sets of empirical DTM are union of around $\binom{n}{q}$ balls with $q=hn$, which makes their computation intractable in practice. To bypass this issue, approximations of the empirical DTM have been proposed in \cite{Merigot2} ($q$-witnessed distance) and \cite{Buchet16} (power distance). Up to our knowledge, these are the only available approximations of the empirical DTM.
The sublevel sets of these two approximations are union of $n$ balls. Thus, it makes the computation of topological invariants more tractable for small data sets, from alpha-shape for instance; see \cite{Ede92}.
Nonetheless, when $n$ is large, there is still a need for a coreset allowing to efficiently compute an approximation of the DTM, as pointed out in \cite{Phillips}. In \cite{MerigotLB}, M\'erigot proves that such a coreset cannot be too small for large dimension.


\subsection{Contribution}

This paper aims at providing such a coreset for the DTM, to face the case where there are many observations, possibly corrupted by noise. We introduce the \textit{$k$-power distance to a measure} $P$ ($k$-PDTM), which is defined as the square root of one of the best $k$-power functions approximating the square of the DTM from above, for the $L_1(P)$ norm. Roughly, we intend to approximate the DTM of a point $x$ with a power distance $\dd_{P,h,k}(x)$ of the form
          \begin{align*}
          \dd_{P,h,k}(x)  = \sqrt{\min_{i\in[\![1,k]\!]} \|x-\theta_i\|^2+\omega^2_{P,h}(\theta_i)},
          \end{align*}
 where the $\theta_i$'s and corresponding $\omega$'s are suitably chosen. Its sub-level sets are union of $k$ balls. Thus, the study of the associated topological invariants gets tractable in practice, even for massive data. 
 

We begin by providing some theoretical guarantees on the $k$-PDTM we introduce. For instance, we prove that it can be expressed as a power distance from a coreset of $k$ points that are local means of the measure $P$. The proofs rely on a geometric study of local sub-measures of $P$ with fixed mass $h\in[0,1]$, showing that such a coreset makes sense whenever $P$ is supported on a compact set. In particular, we prove that the set of means of local sub-measures of $P$ is convex. The discrete case relies on the duality between a weighted Delaunay diagram and its associated weighted Vorono\"i diagram. 

Once the $k$-PDTM properly defined, the main contribution of our paper are the following. First we assess that the $k$-DTM is a good approximation of the DTM in the $L_1$ sense (Proposition \ref{prop:approx_kdtm_dtm}), showing for instance that whenever $M$ has dimension $d'$
\begin{align*}
P \left ( \dd^2_{P,h,k}(u) - \dd^2_{P,h}(u) \right ) \leq C_{P,h} k^{-\frac{2}{d'}}, 
\end{align*}
where $P f(u)$ stands for the integration of $f$ with respect to measure $P$. As mentioned in Proposition \ref{prop:bound_ktdm_dtosupport}, this allows to infer topological guarantees from the sublevel sets of the $k$-PDTM.

Second we prove that this $k$-PDTM shares the robustness properties of the DTM with respect to Wasserstein deformations (Proposition \ref{prop approximation Wasserstein}). Namely, if $Q$ is a sub-Gaussian deformation of $P$ such that the Wasserstein distance $W_2(P,Q) \leq \sigma \leq K$, it holds
\begin{align*}
P \left |\dd^2_{Q,h,k}(u) - \dd_{P,h}^2(u) \right | \leq P \left ( \dd^2_{P,h,k}(u) - \dd^2_{P,h}(u) \right ) + C_{P,h} \sigma K,
\end{align*}
ensuring that the approximation guarantees of our $k$-PDTM are stable with respect to Wasserstein noise. Similar to the DTM, this also guarantees that an empirical $k$-PDTM, that is built on $\X_n$, is a consistent approximation of the true $k$-PDTM. 

At last, we provide more insights on the construction of the empirical $k$-PDTM from a point cloud $\X_n$, facing the practical situation where only a corrupted sample is at hand. We expose a $k$-means like algorithm with complexity $O(n^2hkd)$, and we analyze the approximation performance of such an empirical output. Theorem \ref{thm:prox_empirical_version} shows that, with high probability, 
\begin{align*}
P \left ( \dd^2_{P_n,h,k}(u) - \dd^2_{P,h,k}(u) \right ) \leq C_{P,h} \frac{\sqrt{k}\left(\log(n)\right)^{\frac32}}{\sqrt{n}}.
\end{align*} 
Combining this estimation result with the approximation results between $k$-PDTM and DTM mentioned above suggest that an optimal choice for $k$ is $k=n^{\frac{d'}{d'+4}}$, whenever $M$ has dimension $d'$, resulting in a deviation between empirical $k$-PDTM and DTM of order $n^{-1/(d'+4)}$. This has to be compared with the $n^{-1/d'}$ approximation that the empirical DTM achieves in such cases. In the case where $n$ is large, this $n^{-1/(d'+4)}$ approximation suffices for topological inference. Thus, topological inference built on significantly less points might provide almost similar guarantees than the DTM.

\subsection{Organization of the paper}

This paper is organized as follows. 
In Section 2, we recall some definitions for the DTM that can be expressed as a power distance, and study the set of local means. Section 3 is devoted to the $k$-PDTM, a $k$-power distance which approximates the DTM. We make the link with two equivalent definitions for the $k$-PDTM, derive some stability results, prove its proximity to the DTM highlighting its interest for topological inference. The case of noisy point clouds is addressed in Section 4, where an algorithm to approximate the $k$-PDTM comes up with theoretical guarantees.


\section{Some background about the DTM}\label{section def dtm}
\subsection{Notation and definitions for the DTM}
In the paper, we denote by $\R_d=\{x=(x_1,x_2,\ldots,x_d)\mid\forall i\in[\![1,d]\!],\,x_i\in\R\}$  the $d$-dimensional space equipped with the Euclidean norm $\|.\|$. For $k\in\N^*$ and any space $\A$, $\A^{(k)}$ stands for $\{t=(t_1,t_2,\ldots,t_k)\mid\forall i\in[\![1,k]\!],\,t_i\in\A\}$, where two elements are identified whenever they are equal up to a permutation of the coordinates. Also, $\Sph(0,r)=\{x\in\R_d\mid\|x\|=r\}$ denotes the Euclidean sphere of radius $r$, $\B(x,r)=\{y\in\R_d\mid\|x-y\|<r\}$ the Euclidean ball centred at $x$, and for $c\in\R$ and $v\in\Sph(0,1)$, $\HP(v,c)$ denotes the half-space $\{x\in\R_d\mid\langle x,v\rangle>c\}$. Also, for any subset $A$ of $\R_d$, $\overline{A}$ stands for its closure, $A^\circ$ for its interior, $\partial A=\overline{A}\backslash A^\circ$ its boundary and $A^c=\R_d\backslash A$ its complementary set in $\R_d$. 

In the following, $\PPP(\R_d)$ stands for the set of Borel probability distributions $P$, with support $\supp(P)\subset\R_d$, and, for any $P$-integrable function $f$, $Pf(u)$ denotes the expectation of $f$ with respect to $P$. The following sets of distributions are of particular interest: we denote by $\PPP^K(\R_d)=\left\{P\in\PPP(\R_d)\mid\supp(P)\subset\overline{\B}(0,K)\right\}$ for $K>0$, and $\PPP^{K,h}(\R_d)$ is the set of $P\in\PPP^K(\R_d)$ which put mass neither on the boundaries of balls nor on the half-spaces of $P$-mass $h$. We also allow perturbations of measures in $\PPP^{K,h}(\R_d)$. A sub-Gaussian measure $Q$ with variance $V^2>0$ is a measure $Q\in\PPP(\R_d)$ such that $Q(\B(0,t)^c)\leq\exp(-\frac{t^2}{2V^2})$ for all $t>V$. The set of such measures is denoted by $\PPP^{(V)}(\R_d)$. As well we can define $\PPP^{(V),h}(\R_d)$. The set $\PPP^{(V),h}(\R_d)$ might be thought of as perturbations of $\PPP^{K,h}(\R_d)$. Indeed, if $X=Y+Z$, where $X$ has distribution in $\PPP^{K,h}(\R_d)$ and $Z$ is Gaussian with variance $\sigma^2$, then $Z$ has distribution in $\PPP^{(V),h}(\R_d)$, with $V = K + \sigma$. All these sets of distributions are included in $\PPP_2(\R_d)$, that denotes the set of distributions with finite second moment.


For all $P\in\PPP(\R_d)$ and $n\in\N^*$, $\X_n=\{X_1,X_2,\ldots,X_n\}$ denotes a $n$-sample from $P$, meaning that the $X_i$'s are independent and sampled according to $P$. Also, $P_n=\frac{1}{n}\sum_{i=1}^n\delta_{X_i}$ denotes the empirical measure associated to $P$, where $\delta_{x}\in\PPP(\R_d)$ is such that $\delta_x(\{x\})=1$. Then $\PPP_n(\R_d)$ is the set of $P\in\PPP(\R_d)$ uniform on a set of $n\in\N^*$ points.

An alternative definition to Definition \ref{def:DTM}, for the distance to measure, might be stated in terms of sub-measures. Let $x \in \mathbb{R}^d$. We define $\PPP_{x,h}(P)$ as the set of distributions $P_{x,h}=\frac{1}{h}Q$, for $Q$ a sub-measure of $P$ coinciding with $P$ on $\B(x,\delta_{P,h}(x))$, and such that $Q(\R_d)=h$ and $\supp(Q)\subset\overline{\B}(x,\delta_{P,h}(x))$. Note that when $P\in\PPP^{K,h}(\R_d)$, $\PPP_{x,h}(P)$ is reduced to a singleton $\left\{P_{x,h}\right\}$ with $P_{x,h}$ defined for all Borel sets $B$ by $P_{x,h}(B)=\frac{1}{h}P\left(B\cap\B(x,\delta_{P,h}(x))\right)$. From \cite[Proposition 3.3]{Merigot1}, it holds, for any $x \in \mathbb{R}^d$ and $P_{x,h} \in \PPP_{x,h}(P)$, 
\begin{equation}
\dd_{P,h}^2(x) = P_{x,h}\|x-u\|^2=\|x-m(P_{x,h})\|^2 + v(P_{x,h}),\label{equation dtm Pxh}
\end{equation}
with $m(P_{x,h})=P_{x,h}u$ the mean of $P_{x,h}$ and $v(P_{x,h})=P_{x,h}\|u-m(P_{x,h})\|^2$ its variance.
For convenience, we denote by $M(P_{x,h})=P_{x,h}\|u\|^2$ the second moment of $P_{x,h}$, so that $M(P_{x,h}) = \|m(P_{x,h})\|^2 + v(P_{x,h})$. Whenever $P$ is in $\PPP^{(V)}(\R_d)$, $M$ satisfies the following property.
\lm\label{lm_ m bounded for subgaussian}
Let $P\in\PPP^{(V)}(\R_d)$, then $\forall\ x\in\overline{\R_d}$ and $h\in(0,1]$, $M(P_{x,h})\leq \frac{2V^2}{h}$.
\elm
The proof of Lemma \ref{lm_ m bounded for subgaussian} is deferred to Section \ref{sec:proof_lemma_lm_ m bounded for subgaussian}.

\subsection{From balls to half-spaces: structure of the local means set}
In the previous part, we have seen that the DTM $\dd_{P,h}$ is built from sub-measures of $P$ supported on balls of $P$-mass $h$. Now, by making the center of a ball go to $\infty$ along a direction $v\in\Sph(0,1)$ such that the ball keeps a fixed mass $h$, we obtain a sub-measure of $P$ supported on a half-space, as follows.

For $v\in\Sph(0,1)$, we denote by $v_\infty$ the infinite point associated to the direction $v$. It can be seen as a limit point $\lim_{\lambda\rightarrow+\infty}\lambda v$.
Then, we denote $\overline{\R}_d=\R_d\bigcup\left\{v_\infty\mid v\in\Sph( 0 ,1)\right\}$.
Note that we can equip $\overline{\R}_d$ with the metric $\dR$ defined by $\dR(x,y)=\|\phi(x)-\phi(y)\|$, with $\phi(x)=\frac{x}{\sqrt{1+\|x\|^2}}$ when $x\in\R_d$ and $\phi(v_\infty)=v$ for all $v\in\Sph(0,1)$.
Also, for this metric, a sequence $(x_n)_{n\in\N}$ of $\overline{\R}_d$ converges to $v_\infty$ if and only if $\lim_{n\rightarrow+\infty}\|x_n\|=+\infty$ and $\lim_{n\rightarrow+\infty}\frac{x_n}{\|x_n\|}=v$ with the convention $\frac{w_\infty}{\|w_\infty\|}=w$ for all $w\in\Sph( 0,1)$.\\

Let $v\in\Sph(0,1)$, set $c_{P,h}(v)=\sup\{c\in\R\mid P\left(\left\{x\in\R_d\mid\langle x,v\rangle>c\right\}\right)>h\}$. Then, $\HP(v,c_{P,h}(v))$ corresponds to the largest (for the inclusion order) half-space directed by $v$ with $P$-mass at most $h$, which contains all the $\lambda v$'s for $\lambda$ large enough.
\lm\label{lem_cv_ball_plane}
Let $v\in\Sph(0,1)$ and $P\in\PPP(\R_d)$. Assume that $P(\partial\HP(v,c_{P,h}(v)))=0$.
If $x_n=nv$ for all $n\in\N$, then for $P$-almost all $y\in\R_d$, we have:
\[\lim_{n\rightarrow+\infty}\1_{\B(x_n,\delta_{P,h}(x_n))}(y)=\1_{\HP(v,c_{P,h}(v))}(y).\]
If $(x_n)_{n\in\N}$ is a sequence of $\R_d$ such that $\lim_{n\rightarrow+\infty}\dd_{\overline{\R}_d}(x_n,v_\infty)=0$, then, the result holds up to a subsequence.
\elm
The proof of Lemma \ref{lem_cv_ball_plane} is given in the Appendix, Section \ref{sec:proof Lemma lem_cv_ball_plane}.
For all $P\in\PPP_2(\R_d)$, we can generalize the definition of $\PPP_{x,h}(P)$, $P_{x,h}$, $m(P_{x,h})$, $v(P_{x,h})$ and $M(P_{x,h})$ to the elements $x=v_\infty\in\overline{\R}_d\backslash\R_d$ for all $v\in\Sph(0,1)$.
Note that when $P\in\PPP^{K,h}(\R_d)$, $\PPP_{v_\infty,h}(P)$ is reduced to the singleton $\left\{P_{v_\infty,h}\right\}$ with $P_{v_\infty,h}$ equal to $\frac{1}{h}P(B\cap\HP(v,c_{P,h}(v)))$ for all Borel set $B$. Intuitively, he distributions $P_{v_\infty,h}$ behave like extreme points of $\left\{P_{x,h}\mid x\in\R_d\right\}$. This intuition is formalized by the following Lemma. Denote $\MM_h(P)=\left\{m(P_{x,h})\mid x\in\overline{\R}_d\right\}$.
\lm\label{lm enveloppe convexe}
Let $P\in\PPP_2(\R_d)$, the set $\conv(\MM_h(P))$ is equal to $\bigcap_{v\in\Sph(0,1)}H^c(v,\langle m(P_{v_\infty,h}),v\rangle)$.
\elm
A straightforward consequence of Lemma \ref{lm enveloppe convexe} is the following Lemma \ref{inclusion of conv(M)}.
\lm\label{inclusion of conv(M)}
Let $P\in\PPP_2(\R_d)$, then
\[\forall\, 0<h<h'\leq1,\,\conv\left(\MM_{h'}(P)\right)\subset\conv\left(\MM_{h}(P)\right).\]
\elm
The proofs of Lemmas \ref{lm enveloppe convexe} and \ref{inclusion of conv(M)} are to be found in Section \ref{sec:proof_of_lemma_enveloppe_convexe} and \ref{sec:proof_of_lemma_inclusion_of_conv(M)}. A key property of the local means sets $\MM_h(P)$ is convexity. This will be of particular interest in Section \ref{Two equivalent definitions for the $k$-PDTM}. We begin with the finite-sample case.

\lm\label{M(Pn) convex}
Let $P_n\in\PPP_n(\R_d)$ such that $\supp(P_n)$ is a set of $n$ points in general position, as described in \cite[Section 3.1.4]{Boissonat}, meaning that any subset of $\supp(P_n)$ with size at most $d + 1$ is a set of affinely independent points, set $q\in[\![1,n]\!]$.
Then, the set $\MM_{\frac qn}(P_n)$ is convex.
\elm

\pv[Proof of lemma \ref{M(Pn) convex}]
Let 
\[\hat\MM_h(P_n)=\left\{\bar x=\frac{1}{q}\sum_{p\in\Nn_{q,\X_n}(x)}p\mid x\in\overline{\R}_d,\,\Nn_{q,\X_n}(x)\in\NNN_{q,\X_n}(x)\right\},\] with $\NNN_{q,\X_n}(x)$ the collection of all sets of $q$-nearest neighbors associated to $x$. Note that different $\bar x$ may be associated to the same $x$, and also note that $\hat\MM_h(P_n)\subset\MM_h(P_n)$. Moreover, $\conv(\MM_h(P_n))=\conv(\hat\MM_h(P_n))$ since any $m(P_{n\,x,h})$, for $P_{n\,x,h}\in\PPP_{x,h}(P_n)$, can be expressed as a convex combination of the $\bar x$'s.

Then, $\R_d$ breaks down into a finite number of weighted Vorono\"i cells $\CC_{P_n,h}(\bar x)=\{z\in\R_d\mid\|z-\bar x\|^2+\hat\omega^2(\bar x)\leq\|z-\bar y\|^2+\hat\omega^2(\bar y),\,\forall\bar y\in\hat\MM_h(P_n)\}$, with $\hat\omega^2(\bar x)=\frac{1}{q}\sum_{p\in\Nn_{q,\X_n}(x)}\|p-\bar x\|^2$ the weight associated to any point $\bar x=\frac{1}{q}\sum_{p\in\Nn_{q,\X_n}(x)}p$ in $\hat\MM_h(P_n)$. According to \cite[Theorem 4.3]{Boissonat}, the weighted Delaunay triangulation partitions the convex hull of any finite set of weighted points $\X$ in general position by $d$-dimensional simplices with vertices in $\X$, provided that the associated weighted Vorono\"i cells of all the points in $\X$ are non empty. By duality, (also see \cite[Lemma 4.5]{Boissonat}) these vertices are associated to weighted Vorono\"i cells that have non-empty common intersection.
Thus, any $\theta\in\conv(\MM_h(P_n))$ satisfies $\theta=\sum_{i=0}^{d}\lambda_i\bar x^i$ for some $\bar x^i$'s in $\hat\MM_h(P_n)$ and some non negative $\lambda_i$'s such that $\sum_{i=0}^{d}\lambda_i=1$. 
Also, there exists some $x^*$ in the intersection of the $d+1$ weighted Vorono\"i cells, $\left(\CC_{P_n,h}(\bar x^i)\right)_{i\in[\![0,d]\!]}$.


Set $P_{n\,x^*,h}:=\sum_{i=0}^d\lambda_iP_i$, with $P_i=\frac{1}{q}\sum_{p\in\Nn^i_{q,\X_n}(x^*)}\delta_{\left\{p\right\}}$ when $\bar x^i=\frac{1}{q}\sum_{p\in\Nn^i_{q,\X_n}(x^*)}p$. Then, $P_{n\,x^*,h}$ is a probability measure such that $hP_{n\,x^*,h}$ ($h=\frac{q}{n}$) coincides with $P_n$ on $\B(x,\delta_{P_n,h}(x))$ and is supported on $\overline{\B}(x,\delta_{P_n,h}(x))$. Thus it belongs to $\PPP_{x^*,h}(P_n)$. Moreover, its mean $m(P_{n\,x^*,h})=\theta$. Thus, $\theta\in\MM_h(P_n)$.
\epv
If $P \in \PPP^{K}(\R_d)$, convexity of $\MM_h(P)$ might be deduced from the above Lemma \ref{M(Pn) convex} using the convergence of the empirical distribution $P_n$ towards $P$ in a probabilistic sense. This is summarized by the following Lemma.
\lm\label{lem:transition_M(Pn)_M(P)}
Let $P\in\PPP^K(\R_d)$ and $\theta \in \conv(\MM_h(P))$.
There exists sequences $q_n \in \N$, $\alpha_n \rightarrow 0$,  $P_n\in\PPP_n(\R_d)$ with the points in $\supp(P_n)$ in general position, and $y_n \in \conv(\MM_{\frac{q_n}{n}}(P_n))$ such that
\begin{itemize}
\item[$i)$] $\frac{q_n}{n} \rightarrow h$,
\item[$ii)$] $\|y_n - \theta \| \leq \alpha_n$,
\item[$iii)$] $\sup_{x \in \overline{\R}_d}{\|m(P_{n x,\frac{q_n}{n}})-m(P_{x,h})\|} \leq \alpha_n$. 
\end{itemize} 
\elm        
Lemma \ref{lem:transition_M(Pn)_M(P)} follows from probabilistic arguments when $\X_n$ is sampled at random. Its proof can be found in Section \ref{sec:proof_of_lemma_transition_M(Pn)_M}. Equipped with Lemma \ref{lem:transition_M(Pn)_M(P)}, we can prove the convexity of $\MM_h(P)$.
\prop\label{M(P) convex}
If $P\in\PPP^{K}(\R_d)$ for $K>0$ is such that $P(\partial\HP(v,c))=0$ and $P(\partial\B(x,r))=0$ for all $v\in\Sph(0,1)$, $x\in\R_d$, $c\in\R$, $r\geq0$, then for all $h\in(0,1]$,
$\MM_{h}(P)$ is convex.
\eprop
\pv[Proof of Proposition \ref{M(P) convex}]
Let $\theta\in\conv(\MM_h(P))$, $P_n$, $q_n$, $\alpha_n$, $y_n$ as in Lemma \ref{lem:transition_M(Pn)_M(P)} and for short let $h_n = \frac{q_n}{n}$.



Since $\MM(P_n,h_n)$ is convex, there is a sequence $(x_n)_{n\geq N}$ in $\overline{\R}_d$ such that $y_n=m(P_{n\,x_n,h_n})$ converges to $\theta$. If $(x_n)_{n\geq N}$ is bounded, then up to a subsequence we have $x_n \rightarrow x$, for some $x \in \mathbb{R}^d$. If not, considering $\frac{x_n}{\|x_n\|}$, up to a subsequence we get $x_n \rightarrow v_\infty$. In any case $x_n \rightarrow x$, for $x \in \overline{\R}_d$. Combining  Lemma \ref{lem:transition_M(Pn)_M(P)} and Lemma \ref{lem_cv_ball_plane} yields $\theta=m(P_{x,h})$. Thus, $\theta\in\MM_h(P)$.
\epv 


\subsection{The DTM defined as a power distance}

A \textbf{power distance} indexed on a set $I$ is the square root of a power function $f_{\tau,\omega}$ defined on $\R_d$ from a family of centers $\tau=(\tau_i)_{i\in I}$ and weights $\omega=(\omega_i)_{i\in I}$ by $f_{\tau,\omega}:x\mapsto\inf_{i\in I}\|x-\tau_i\|^2+{\omega_i}^2$. A \textbf{$k$-power distance} is a power distance indexed on a finite set of cardinal $|I|=k$.

In \cite[Proposition 3.3]{Merigot1}, the authors point out that $P_{x,h}\|x-u\|^2\leq Q\|x-u\|^2$ for all $Q\in\PPP(\R_d)$ such that $hQ$ is a sub-measure of $P$. 
This remark, together with \eqref{equation dtm Pxh}, provides an expression for the DTM as a power distance.
\prop[{\cite[Proposition 3.3]{Merigot1}}]\label{DTM as power function}
If $P\in\PPP_2(\R_d)$, then for all $x\in\R_d$, we have:
\[\dd_{P,h}^2(x)=\inf_{y\in\overline{\R}_d}\inf_{P_{y,h}\in\PPP_{y,h}(P)}\|x-m(P_{y,h})\|^2+v(P_{y,h}),\]
and the infimum is attained at $y=x$ and any measure $P_{x,h}\in\PPP_{x,h}(P)$.
\eprop

As noted in M\'erigot et al \cite{Merigot2}, this expression holds for the empirical DTM $\dd_{P_n,h}$. In this case, $m(P_{n,x,h})$ corresponds to the barycentre of the $q=nh$ nearest-neighbors of $x$ in $\X_n$, $\Nn_{q,\X_n}(x)$, and $v(P_{n,y,h})=\frac{1}{q}\sum_{p\in\Nn_{q,\X_n}(x)}\|x-p\|^2$, at least for points $x$ whose set of $q$ nearest neighbors is uniquely defined.

\subsection{Semiconcavity and DTM}

In the following, we will often use the following lemma connected to the property of concavity of the function $x\mapsto \dd^2_{P,h}(x)-\|x\|^2$.
\lm[{\cite[Proposition 3.6]{Merigot1}}]
\label{lemme de concavite}
If $P\in\PPP(\R_d)$, then for all $x,\,y\in\R_d$ and $P_{x,h}\in\PPP_{x,h}(P)$,
\[\dd^2_{P,h}(y)-\|y\|^2\leq\dd^2_{P,h}(x)-\|x\|^2-2\langle y-x,m(P_{x,h})\rangle,\]
with equality if and only if $P_{x,h}\in\PPP_{y,h}(P)$.\\
\elm

\section{The $k$-PDTM: a coreset for the DTM}\label{sec:kPDTM_a_coreset_for_DTM}

In Proposition \ref{DTM as power function}, we have written the DTM as a power distance. This remark has already been exploited in \cite{Merigot2} and \cite{Buchet16}, where the DTM has been approximated by $n$-power distances.
In this paper, we propose to keep only $k$ centers.
\df\label{def:Opt}
For any $P\in\PPP_2(\R_d)$, we define $\opt(P,h,k)$ by: 
\[\opt(P,h,k)=\argmin\left\{P\min_{i\in[\![1,k]\!]}\|u-m(P_{t_i,h})\|^2+v(P_{t_i,h})\mid t=(t_1,t_2,\ldots t_k)\in\overline{\R}_d^{(k)}\right\}.\]
\edf
A closely related notion to Definition \ref{def:Opt} is the following weighted Voronoï measures.
\df
A set of \textbf{weighted Vorono\"i measures} associated to a distribution $P\in\PPP_2(\R_d)$, $t\in\overline{\R}_d^{(k)}$ and $h\in(0,1]$ is a set $\left\{\tilde P_{t_1,h}, \tilde P_{t_2,h},\ldots\tilde P_{t_k,h}\right\}$ of $k\in\N^*$ positive sub-measures of $P$ such that $\sum_{i=1}^k\tilde P_{t_i}=P$ and 
\[\forall x\in\supp(\tilde P_{t_i,h}),\,\|x-m(P_{t_i,h})\|^2+v(P_{t_i,h})\leq\|x-m(P_{t_j,h})\|^2+v(P_{t_j,h}),\,\forall j\in[\![1,k]\!].\]
We denote by $\tilde m(\tilde P_{t_i,h})=\frac{\tilde P_{t_i,h}u}{\tilde P_{t_i,h}(\R_d)}$ the expectation of $\tilde P_{t_i,h}$, with the convention $\tilde m(\tilde P_{t_i,h})=0$ when $\tilde P_{t_i,h}(\R_d)=0$.
\edf
Note that a set of weighted Vorono\"i measures can always be assigned to any $P\in\PPP_2(\R_d)$ and $t\in\overline{\R}_d^{(k)}$, it suffices to split $\R_d$ in weighted Vorono\"i cells associated to the centers $(m(P_{t_i,h}))_{i\in[\![1,k]\!]}$ and weights $(v(P_{t_i,h}))_{i\in[\![1,k]\!]}$, see \cite[Section 4.4.2]{Boissonat}, and split the remaining mass on the border of the cells in a measurable arbitrary way.

\thm\label{opt nonempty}
For all $h\in(0,1]$, $k\in\N^*$ and $P\in\PPP^{K}(\R_d)$ for some $K>0$, such that $P(\partial\HP(v,c_{P,h}(v)))=0$ for all $v\in\Sph(0,1)$, the set $\opt(P,h,k)$ is not empty.
Moreover, there is some $s\in\opt(P,h,k)\cap\overline{\B}(0,K)^{(k)}$ such that $s_i=\tilde m(\tilde P_{s_i,h})$ for all $i\in[\![1,k]\!]$.
\ethm
\pv[Sketch of proof]
For $s\in\overline{\R}_d^{(k)}$, set $f_s:x\in\R_d\mapsto\min_{i\in[\![1,k]\!]}\left(\|x-m(P_{s_i,h})\|^2+v(P_{s_i,h})\right)$. Then, Lemma \ref{lem_cv_ball_plane} and the dominated convergence theorem yield $\inf_{t\in\R_d}Pf_t(u)=\inf_{t\in\overline{\R}_d}Pf_t(u)$.

Let $(t_n)_{n\in\N}$ be a sequence in $\R_d^{(k)}$ such that $Pf_{t_n}(u)\leq\inf_{t\in\overline{\R}_d}Pf_t(u)+\frac{1}{n}$, and denote by $m^*$ the limit of a converging subsequence of $\left(\tilde m(\tilde P_{t_{n,1},h}),\tilde m(\tilde P_{t_{n,2},h}),\ldots,\tilde m(\tilde P_{t_{n,k},h})\right)_{n\in\N}$ in the compact space $\overline{\B}(0,K)^{(k)}$.
Then, Lemma \ref{lemme de concavite} and (\ref{equation dtm Pxh}) yield $Pf_{m^*}(u)=\inf_{t\in\overline{\R}_d}Pf_t(u)$.

Set $s_i=\tilde m(\tilde P_{m^*_{i},h})$ for all $i\in[\![1,k]\!]$, then $\tilde m(\tilde P_{s_i,h})=s_i$ and $Pf_{m^*}(u)=Pf_s(u)$.
\epv
The detailed proof of Theorem \ref{opt nonempty} is given in Section \ref{sec:proof_thm_opt_nonempty}. Note that the distributions in $\PPP^{K,h}$ are in the scope of Theorem \ref{opt nonempty}. 
\subsection{Two equivalent definitions for the $k$-PDTM}
\label{Two equivalent definitions for the $k$-PDTM}

\df
Let $P\in\PPP_2(\R_d)$, the \textbf{$k$-power distance to a measure} ($k$-PDTM) $\dd_{P,h,k}$ is defined for any $s\in\opt(P,h,k)$ by:
\[\dd^2_{P,h,k}(x)=\min_{i\in[\![1,k]\!]} \|x-m(P_{s_i,h})\|^2+v(P_{s_i,h}).\]
An \textbf{$\epsilon$-approximation of the $k$-PDTM}, denoted by an $\dd^2_{P,h,k,\epsilon}$ is a function defined by the previous expression but for some $s\in\overline{\R}_d^{(k)}$ satisfying
\[P\min_{i\in[\![1,k]\!]} \|u-m(P_{s_i,h})\|^2+v(P_{s_i,h})\leq\inf_{t\in\overline{\R}_d^{(k)}}P\min_{i\in[\![1,k]\!]}\|u-m(P_{t_i,h})\|^2+v(P_{t_i,h})+\epsilon.\]
\edf

Theorem \ref{opt nonempty} states that the $k$-PDTM is well defined when $P\in\PPP^K(\R_d)$ and satisfies $P(\partial\HP_{P,h}(v,c_{P,h}(v)))=0$ for all $v\in\Sph(0,1)$. Nonetheless, whenever $\opt(P,h,k)$ is not a singleton, the $k$-PDTM is not unique.
Note that for all $x\in\R_d$, $\dd_{P,h,k}(x)\geq\dd_{P,h}(x)$.
\df
The set $\OPT(P,h,k)$ is defined by:
\[\OPT(P,h,k)=\argmin\left\{P\min_{i\in[\![1,k]\!]}\|u-\tau_i\|^2+\omega^2_{P,h}(\tau_i)\mid\tau=(\tau_1,\tau_2,\ldots \tau_k)\in\R_d^{(k)}\right\},\]
with $\omega_{P,h}(\tau)=\inf\left\{\omega>0\mid\forall x\in\R_d,\,\|x-\tau\|^2+\omega^2\geq\dd^2_{P,h}(x)\right\}$ for $\tau\in\R_d$, that is:
\begin{equation}
\omega^2_{P,h}(\tau)=\sup_{x\in\R_d}\dd^2_{P,h}(x)-\|x-\tau\|^2.\label{equation omega square}
\end{equation}
\edf
The following Lemma shows that $\OPT(P,h,k)$ is included in $\conv(\MM_h(P))^{(k)}$.
\lm\label{lem:omega_convM}
Let $P\in\PPP^K(\R_d)\cup\PPP^{(V)}(\R_d)$. Then $\theta\in\conv(\MM_h(P))$ if and only if $\omega_{P,h}(\theta) <+\infty$.
\elm
\pv[Proof of Lemma \ref{lem:omega_convM}]
According to Proposition \ref{DTM as power function}, for all $x\in\R_d$, $\dd^2_{P,h}(x)-\|x-\theta\|^2$ may be written as
\[\inf_{y\in\overline{\R}_d}\inf_{P_{y,h}\in\PPP_{y,h}(P)}\left\{\|m(P_{y,h})\|^2 + v(P_{y,h}) - \|\theta\|^2 + 2\langle x,\theta-m(P_{y,h})\rangle\right\},\]
which is lower-bounded by
\[\inf_{y\in\overline{\R}_d}\inf_{P_{y,h}\in\PPP_{y,h}(P)}\left\{\|m(P_{y,h})\|^2+v(P_{y,h})\right\}-\|\theta\|^2+\inf_{\tau\in\MM_h(P)}\left\{2\langle x,\theta-\tau\rangle\right\}.\]
Assume $\theta\notin\conv(\MM_h(P))$. According to Lemma \ref{lm enveloppe convexe}, $\conv(\MM_h(P))$ is a convex and compact subset of $\R_d$.
The Hahn-Banach separation theorem thus provides some vector $v\in\R_d$ and $C>0$ such that $\forall\tau\in\MM_h(P)$, $\langle\theta-\tau,v\rangle<C$. Setting $x_n=-nv$ for $n\in\N^*$ yields $\lim_{n\rightarrow+\infty}\inf_{\tau\in\MM_h(P)}\langle x_n,\theta-\tau\rangle=+\infty$. Thus, $\sup_{x\in\R_d}\dd^2_{P,h}(x)-\|x-\theta\|^2=+\infty$.

Now, let $\theta\in\conv(\MM_h(P))$, we can write $\theta=\sum_{i=0}^d\lambda_i m(P_i)$ for $P_i=P_{x_i,h}$ with the $x_i$'s in $\overline{\R}_d$. We have:
\begin{align*}
\sup_{x\in\overline{\R}_d}\dd^2_{P,h}(x)-\|x-\theta\|^2 & = \sup_{x\in\overline{\R}_d}\sum_{i=0}^{d}\lambda_i(\dd^2_{P,h}(x)-\|x-\theta\|^2)\\
 &\leq\sup_{x\in\overline{\R}_d}\sum_{i=0}^d\lambda_i\left(\|x-m(P_i)\|^2+v(P_i)-\|x-\theta\|^2\right)\\
 &=\sup_{x\in\overline{\R}_d}\sum_{i=0}^d\lambda_i\left(v(P_i)+2\langle x,\theta-m(P_i)\rangle+\|m(P_i)\|^2-\|\theta\|^2\right)\\
 &=\sum_{i=0}^d\lambda_i\left(v(P_i)+\|m(P_i)\|^2-\|\theta\|^2\right),
\end{align*}
according to Proposition \ref{DTM as power function}.
Thus, we get that
\begin{align}\label{eq:bound_for_omega}
\omega^2_{P,h}(\theta)+\|\theta\|^2&\leq\sum_{i=0}^d\lambda_i(v(P_i)+\|m(P_i)\|^2)\leq\sup_{x\in\overline{\R}_d}\left\{v(P_{x,h})+\|m(P_{x,h})\|^2\right\}.
\end{align}
Lemma \ref{lm_ m bounded for subgaussian} yields $\omega^2_{P,h}(\theta) < + \infty$.
\epv
\thm\label{th_ equivalence}
If $P\in\PPP^{K,h}(\R_d)$ for some $h\in(0,1]$ and $K>0$, or $P\in\PPP_n(\R_d)$ such that $\supp(P_n)$ is a set of $n$ points in general position as described in Lemma \ref{M(Pn) convex}, for some $h=\frac{q}{n}$ with $q\in[\![1,n]\!]$, then,
any function $\dd_{P,h,k}$ satisfies for some $\theta\in\OPT(P,h,k)$:
\[\dd^2_{P,h,k}(x)=\min_{i\in[\![1,k]\!]} \|x-\theta_i\|^2+\omega^2_{P,h}(\theta_i),\,\forall x\in\R_d.\]
Conversely, for all $\theta\in\OPT(P,h,k)$, $x\mapsto\sqrt{\min_{i\in[\![1,k]\!]} \|x-\theta_i\|^2+\omega^2_{P,h}(\theta_i)}$ is a $k$-PDTM.
\ethm

\pv[Proof of Theorem \ref{th_ equivalence}]
For all $\tau\in\R_d^{(k)}$, for all $i\in[\![1,k]\!]$, if $\tau_i\notin\conv(\MM_h(P))$, then according to Lemma \ref{lem:omega_convM}, $\omega_{P,h}(\tau_i)=+\infty$. In this case, $\tau\notin\OPT(P,h,k)$.\\
Thus, for all $\tau\in\OPT(P,h,k)$, for all $i$, $\tau_i\in\conv(\MM_h(P))$. According to Proposition \ref{M(P) convex} and Lemma \ref{M(Pn) convex}, $\MM_h(P)$ is convex.
Thus, \[\OPT(P,h,k)=\argmin\left\{P\min_{i\in[\![1,k]\!]}\|u-\tau_i\|^2+\omega^2_{P,h}(\tau_i)\mid\tau=(\tau_1,\tau_2,\ldots \tau_k)\in\MM_h(P)^{(k)}\right\}.\]
Moreover, according to Proposition \ref{DTM as power function}, and (\ref{equation omega square}), $\omega_{P,h}^2(m(P_{t,h}))=v(P_{t,h})$, for all $t\in\overline{\R}_d$.
Thus,\[\inf_{t\in\overline{\R}_d^{(k)}}P\min_{i\in[\![1,k]\!]}\|x-m(P_{t_i,h})\|^2+v(P_{t_i,h})=\inf_{\tau\in\R_d^{(k)}}P\min_{i\in[\![1,k]\!]}\|u-\tau_i\|^2+\omega^2_{P,h}(\tau_i).\]
\epv

Therefore, Theorem \ref{th_ equivalence} allows to consider the function $\dd_{P,h,k}$ as the square root of a minimizer of the $L_1(P)$ norm $f\mapsto P|f-\dd^2_{P,h}|(u)$ among all the $k$-power functions $f$ which graph lies above the graph of the function $\dd^2_{P,h}$.\\

\subsection{Proximity to the DTM}

Here we show that the $k$-PDTM approximates the DTM in the following sense.

\prop\label{prop:approx_kdtm_dtm}
Let $P\in\PPP^K(\R_d)$ for $K>0$ and let $M \subset \B(0,K)$ be such that $P(M)=1$. Let $f_M(\varepsilon)$ denote the $\varepsilon$ covering number of $M$. Then we have
\[0\leq P\dd^2_{P,h,k}(u)-\dd^2_{P,h}(u)\leq 2 f_{M}^{-1}(k) \zeta_{P,h}(f_{M}^{-1}(k)), \quad \mbox{with} \quad f^{-1}_M(k) = \inf \left \{ \varepsilon >0| \quad f_{M}(\varepsilon) \leq k \right \},\]
where $\zeta_{P,h}$ is the continuity modulus of $x\mapsto m(P_{x,h})$, that is
\[\zeta_{P,h}(\epsilon)= \sup_{x,y \in M, \|x-y\| \leq \varepsilon}\left\{|m(P_{x,h})-m(P_{y,h})|\right \}.
\]
\eprop

\pv[Proof of Proposition \ref{prop:approx_kdtm_dtm}]
The first inequality comes from Proposition \ref{DTM as power function}.

We then focus on the second bound.
By definition of $\dd_{P,h,k}$, for all $x\in\R_d$ and $t=(t_1,t_2,\ldots,t_k)\in\R_d^{(k)}$ we have:
$P\dd^2_{P,h,k}(x)\leq P\min_{i\in[\![1,k]\!]} \|u-m(P_{t_i,h})\|^2+v(P_{t_i,h}).$
Thus,
\begin{align*}
P\dd^2_{P,h,k}(u)-\dd^2_{P,h}(u)&\leq P\min_{i\in[\![1,k]\!]}\|u-m(P_{t_i,h})\|^2+v(P_{t_i,h})-\dd^2_{P,h}(u)\\
&=P\min_{i\in[\![1,k]\!]}(\dd^2_{P,h}(t_i)-\|t_i\|^2)-(\dd^2_{P,h}(u)-\|u\|^2)+\langle u-t_i,-2m(P_{t_i,h})\rangle\\
&\leq P\min_{i\in[\![1,k]\!]}2\langle u-t_i,m(P_{u,h})-m(P_{t_i,h})\rangle \\
& \leq 2 P \min_{i\in[\![1,k]\!]} \|u-t_i\| \|m(P_{u,h})-m(P_{t_i,h})\|,
\end{align*}
where we used (\ref{equation dtm Pxh}), Lemma \ref{lemme de concavite} and Cauchy-Schwarz inequality.
Now choose $t_1, \hdots, t_k$ as a $f_{M}^{-1}(k)$-covering of $M$. The result follows.      
\epv

When $P$ is roughly uniform on its support, the quantities $f_M^{-1}(k)$ and $\zeta_{P,h}$ mostly depend on the dimension and radius of $M$. We focus on two cases in which Proposition \ref{prop:approx_kdtm_dtm} may be adapted. First, the case where the distribution $P$ has an ambient-dimensional support is investigated.
\cor\label{cor:approx_lebesgue}
   Assume that $P$ have a density $f$ satisfying $0<f_{min} \leq f \leq f_{max}$. Then
   \[
   0\leq P\dd^2_{P,h,k}(u)-\dd^2_{P,h}(u) \leq C_{f_{max},K,d,h} k^{-2/d}.
   \]
\ecor
The proof of Corollary \ref{cor:approx_lebesgue} is given in Section \ref{sec:proof_cor_approx_lebesgue}. Note that no assumptions on the geometric regularity of $M$ is required for Corollary \ref{cor:approx_lebesgue} to hold. In the case where $M$ has a lower-dimensional structure, more regularity is required, as stated by the following corollary.

\cor\label{cor:approx_manifold}
       Suppose that $P$ is supported on a compact $d'$-dimensional $\mathcal{C}^2$-submanifold of $\B(0,K)$, denoted by $N$.  Assume that $N$ has positive reach $\rho$, and that $P$ has a density $0 < f_{min} \leq f \leq f_{max}$ with respect to the volume measure on $N$.      
       Moreover, suppose that $P$ satisfies, for all $x \in N$ and positive $r$, 
       \begin{equation}\label{eq:abstandard1}
       P(\B(x,r)) \geq c f_{min} r^{d'}\wedge 1. 
       \end{equation}       
       Then, for $k\geq c_{N,f_{min}}$ and $h \leq c_{N,f_{min}}$, we have $ 0\leq P\dd^2_{P,h,k}(u)-\dd^2_{P,h}(u) \leq C_{N,f_{min},f_{max}} k^{-2/d'}.$
\ecor

Note that (\ref{eq:abstandard1}), also known as $(cf_{min},d')$-standard assumption, is usual in set estimation (see, e.g., \cite{Chazal15}). In the submanifold case, it may be thought of as a condition preventing the boundary from being arbitrarily narrow. This assumption is satisfied for instance in the case where $\partial N$ is empty or is a $\mathcal{C}^2$ $d'-1$-dimensional submanifold (see, e.g., \cite[Corollary 1]{Aaron16}). An important feature of Corollary \ref{cor:approx_manifold} is that this approximation bound does not depend on the ambient dimension. The proof of Corollary \ref{cor:approx_manifold} may be found in Section \ref{sec:proof_cor_approx_manifold}.

\subsection{Wasserstein stability for the $k$-PDTM}
Next we assess that our $k$-PDTM shares with the DTM the key property of robustness to noise.
\prop\label{prop approximation Wasserstein}
Let $P\in\PPP^K(\R_d)$ for some $K>0$, $Q\in\PPP_2(\R_d)$, and $\epsilon>0$. Set $\dd^2_{Q,h,k,\epsilon}$ an $\epsilon$-approximation of the $k$-PDTM of $Q$, then $P\left|\dd^2_{Q,h,k,\epsilon}(u)-\dd^2_{P,h}(u)\right|$ is bounded from above by $B_{P,Q,h,k,\epsilon}$ with \[B_{P,Q,h,k,\epsilon}=\epsilon + 3\|\dd^2_{Q,h}-\dd^2_{P,h}\|_{\infty,\supp(P)} + P \dd^2_{P,h,k}(u)-\dd^2_{P,h}(u) + 2W_1(P,Q)\sup_{s\in\overline{\R}_d}\|m(P_{s,h})\|.\]
\eprop
Note that Lemma \ref{lm_ m bounded for subgaussian} gives a bound on $m(Q_{s,h})$ whenever $Q$ is sub-Gaussian. Also, upper-bounds for the deviation of the $k$-PDTM to the DTM associated to $P$ have been derived in the previous subsection.

\pv[Sketch of proof]
For $x\in\supp(P)$, $\max\left\{0,-\left(\dd^2_{Q,h,k}-\dd^2_{P,h}\right)(x)\right\}\leq\|\dd^2_{P,h}-\dd^2_{Q,h}\|_{\infty,\supp(P)}$.
Set $f_{Q,p}(x)=2\langle x,m(Q_{p,h})\rangle+v(Q_{p,h})$ for $p\in\R_d$, and let $t\in\opt(Q,h,k)$. Then, $P-Q \min_{i\in[\![1,k]\!]}f_{Q,t_i}(u)\leq 2W_1(P,Q)\sup_{t\in\overline{\R}_d}m(Q_{t,h})$ and for $s$ given by Theorem \ref{opt nonempty}, that is $s\in\opt(P,h,k)\cap\overline{\B}(0,K)^{(k)}$ such that $s_i=\tilde m(\tilde P_{s_i,h})$ for all $i\in[\![1,k]\!]$, $P \min_{i\in[\![1,k]\!]}f_{Q,s_i}(u)-\min_{i\in[\![1,k]\!]}f_{P,s_i}(u)$ is bounded from above by $\|\dd^2_{P,h}-\dd^2_{Q,h}\|_{\infty,\supp(P)}$ .
\epv
The details of the proof of Proposition \ref{prop approximation Wasserstein} can be found in Section \ref{sec:proof_prop_approximation_wasserstein}.


\subsection{Geometric inference with the $k$-PDTM}\label{sec:geominf_kdtm} 

As detailed in \cite[Section 4]{Merigot1}, under suitable assumptions, the sublevel sets of the distance to measure are close enough to the sublevel sets of  the distance to its support. Thus they allow to infer the geometric structure of the support. As stated below, this is also the case when replacing the distance to measure with the $k$-PDTM.

\prop\label{prop:bound_ktdm_dtosupport}
Let $M$ be a compact set in $\B(0,K)$ such that $P(M)=1$. Moreover, assume that there exists $d'$ such that, for every $p \in M$ and $r\geq 0$, 
\begin{align}\label{eq:abstandard}
P(\B(p,r)) \geq C(P)r^{d'} \wedge 1.
\end{align}
Let $Q$ be a probability measure (thought of as a perturbation of $P$), and let $\Delta_P^2$ denote $P\dd^2_{Q,h,k,\varepsilon}(u)$.
Then, we have 
\[
\sup_{x \in \R_d} |d_{Q,h,k,\varepsilon}(x) - d_M(x)| \leq  C(P)^{-\frac{1}{d'+2}} \Delta_P^{\frac{2}{d'+2}} +  W_2(P,Q)h^{-\frac{1}{2}},
\]
where $W_2$ denotes the Wasserstein distance.
\eprop
Proposition \ref{prop:bound_ktdm_dtosupport}, whose proof can be found in Section \ref{sec:proof_bound_fdtm_dtosupport}, ensures that the $k$-PDTM achieves roughly the same performance as the distance to measure (see, e.g., \cite[Corollary 4.8]{Merigot1}) provided that $d^2_{Q,h,k,\varepsilon}$ is small enough on the support $M$ to be inferred. As will be shown in the following Section, this will be the case if $Q$ is an empirical measure drawn close to the targeted support.

\section{Approximation of the $k$-PDTM from point clouds}\label{sec:Approximation_kdtm_pointclouds}

Let $P\in\PPP_2(\R_d)$, an approximation of the $k$-PDTM $\dd_{P,h,k}$, is given by the \textbf{empirical $k$-PDTM} $\dd_{P_n,h,k}$.
Note that when $k=n$, $\dd_{P_n,h,n}$ is equal to the $q$-witnessed distance.
Also, when $h=0$ we recover the $k$-means method. 
\subsection{An algorithm for the empirical $k$-PDTM}

The following algorithm is inspired by the Lloyds algorithm.
We assume that the mass parameter $h=\frac{q}{n}$ for some positive integer $q$. And for any $t\in\R_d$, we use the notation $c(t)=\frac{1}{q}\sum_{i=1}^qX_i(t)$, where $X_i(t)$ is one of the $i$-th nearest neighbor of $t$ in $\R_d$. We denote $\omega^2(t)=\frac{1}{q}\sum_{i=1}^{q}\left(X_i(t)-c(t)\right)^2$, and $\CC(t)$ the weighted Vorono\"i cell associated to $t$. We use the notation $|\CC(t)|$ for the cardinal of $\CC(t)\cap\X_n$.

\begin{lstlisting}[frame=single,caption={Local minimum algorithm},label=list:8-6,abovecaptionskip=-\medskipamount]
Input:%*$\X_n$ a $n$-sample from $P$, $q$ and $k$ *);
# Initialization
%* Sample $t_1$, $t_2$,\ldots $t_k$ from $\X_n$ without replacement.*);
while the %*$t_i$*)s vary make the following two steps:
  # Decomposition in weighted Voronoi cells.
  for j in 1..%*$n$*):
    %*Add $X_j$ to the $\CC(t_i)$ (for $i$ as small as possible) satisfying *)
    %*$\|X_j-c(t_i)\|^2+\omega^2(t_i)\leq\|X_j-c(t_l)\|^2+\omega^2(t_l)\,\forall l\neq i$*);
  # Computation of the new centers and weights.
  for i in 1..%*$k$*):
    %*$t_i={{1}\over{|\CC(t_i)|}}\sum_{X\in\CC(t_i)}X$*);
Output:%*$(t_1,t_2,\ldots,t_k)$*)
\end{lstlisting}

The following proposition relies on the same arguments as in the proof of Theorem \ref{opt nonempty}.
\prop\label{prop algo}
This algorithm converges to a local minimum of $t\mapsto P_n\min_{i\in[\![1,k]\!]}\|x-m(P_{n\,t_i,h})\|^2+v(P_{n\,t_i,h})$.
\eprop
The proof of Proposition \ref{prop algo} can be found in Section \ref{sec:proof_prop_algo}.
\subsection{Proximity between the $k$-PDTM and its empirical version}

\thm\label{thm:prox_empirical_version}
Let $P$ be supported on $M \subset \B(0,K)$. Assume that we observe $X_1, \hdots, X_n$ such that $X_i = Y_i + Z_i$, where $Y_i$ is an i.i.d $n$-sample from $P$ and $Z_i$ is sub-Gaussian with variance $\sigma^2$, with $\sigma \leq K$. 
Let $Q_n$ denote the empirical distribution associated with the $X_i$'s.
Then, for any $p>0$, with probability larger than $1-7n^{-p}$, we have
\begin{align*}
 P(\dd^2_{Q_n,h,k}(u) - d^2_{Q,h,k}(u))  \leq C \sqrt{k} \frac{K^2((p+1)\log(n))^{\frac{3}{2}}}{h \sqrt{n}} + C \frac{K \sigma}{\sqrt{h}}.
\end{align*} 
\ethm

A proof of Theorem \ref{thm:prox_empirical_version} is given in Section \ref{sec:proof_thm_prox_empirical_version}. Theorem \ref{thm:prox_empirical_version}, combined with Proposition \ref{prop approximation Wasserstein}, allows to choose $k$ in order to minimize $P\dd^2_{Q_n,h,k}(u)$. Indeed, in the framework of Corollaries \ref{cor:approx_lebesgue} and \ref{cor:approx_manifold} where the support has intrinsic dimension $d'$, such a minimization boils down to optimize a quantity of the form
\[
\frac{C \sqrt{k} K^2((p+1)\log(n))^{\frac{3}{2}}}{h\sqrt{n}} + C_P k^{-\frac{2}{d'}}.
\]
Hence the choice $k \sim n^{\frac{d'}{d'+4}}$ ensures that for $n$ large enough, only $n^{\frac{d'}{d'+4}}$ points are sufficient to approximate well the sub-level sets of the distance to support. For surface inference ($d'=2$), this amounts to compute the distance to $n^{\frac{1}{3}}$ points rather than $n$, which might save some time. Note that when $d'$ is large, smaller choices of $k$, though suboptimal for our bounds, would nonetheless give the right topology for large $n$'s. In some sense, Theorem \ref{thm:prox_empirical_version} advocates only an upper bound on $k$, above which no increase of precision can be expected.

\subsection{Some numerical illustration}
As in \cite{Merigot2}, we sampled $n=6000$ points from the uniform measure on a sideways with radius $\sqrt{2}$ and $\sqrt{{9}\over{8}}$ convolved with a Gaussian $\NN(0,\sigma^2)$ with $\sigma=0.45$. We then plotted in grey the $r$-sub-level set of the $q$-witnessed distance and in purple, the $r$-sub-level set of an approximation of $\dd_{P_n,q,k}$ with $r=0.24$ and $q=50$ nearest-neighbors. The approximation of $\dd_{P_n,q,k}$ is obtained after running our algorithm 10 times and keeping the best (for the $L_1(P_n)$ loss) function obtained, each time after at most 10 iterations.
\begin{figure}[h!]
 \begin{minipage}[h]{.3\linewidth}
  \centering\includegraphics[scale=0.35,bb=150 -10 150 150]{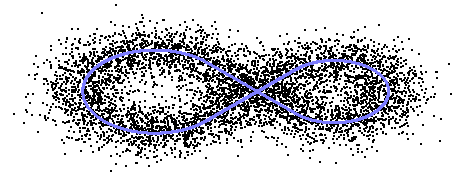}
  \caption{6000-sample}
 \end{minipage} \hfill
 \begin{minipage}[h]{.3\linewidth}
  \centering\includegraphics[scale=0.35,bb=150 -10 150 150]{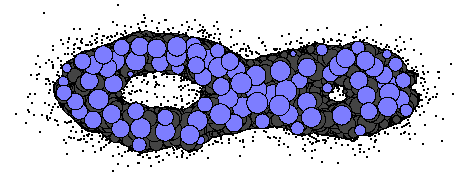}
  \caption{$k=100$}
 \end{minipage} \hfill
 \begin{minipage}[h]{.3\linewidth}
  \centering\includegraphics[scale=0.35,bb=150 -10 150 150]{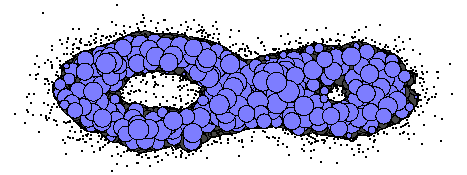}
  \caption{$k=300$}
 \end{minipage}
\end{figure}
Choosing $k=100$ points leads to a too sparse approximation of the sublevel sets of the $q$-witnessed distance. On the contrary, small holes which appeared in the $r$-sub-level set, when $k=300$, will disappear quickly when the radius $r$ will get larger, before the two holes get filled.\\

The authors are grateful to Pascal Massart, Fr\'ed\'eric Chazal and Marc Glisse for their precious advice.

\printbibliography

\pagebreak
\part*{Appendix}
\begin{appendix}

\section{Proofs for Section \ref{section def dtm}}
\subsection{Proof of Lemma \ref{lm_ m bounded for subgaussian}}\label{sec:proof_lemma_lm_ m bounded for subgaussian}
\pv
       Let $P\in\PPP^{(V)}(\R_d)$, $x$ in $\overline{\R}_d$ and $P_{x,h}$ a sub-measure of $P$, supported on $\overline{\B}(x,\delta_{P,h}(x))$ (or on $\HP(v,c_{P,h}(v))$ if $x=v_\infty \in \overline{\R}_d\backslash\R_d$),  coinciding with $P$ on $\B(x,\delta_{P,h}(x))$, and such that $P_{x,h}(\R_d)=h$. We may write
       \begin{align*}
       P_{x,h} \|u\|^2 & \leq \frac{1}{h} P \|u\|^2 \\
                       & \leq \frac{1}{h} \left [P \|u\|^2 \1_{\|u\|\leq V} + P \|u\|^2 \1_{\|u\|> V} \right ] \\
                       & \leq \frac{V^2}{h} + \frac{P \|u\|^2 \1_{\|u\|> V}}{h}.
\end{align*}        
Since $ P \|u\|^2 \1_{\|u\|> V} \leq N_{V^2}t^2 \1_{t> V} \leq   N_{V^2}t^2 = V^2$, where $N_{V^2}$ denotes the distribution of a Gaussian distribution with variance $V^2$, the result of Lemma \ref{lm_ m bounded for subgaussian} follows.   
\epv
\subsection{Proof of Lemma \ref{lem_cv_ball_plane}}\label{sec:proof Lemma lem_cv_ball_plane}

\pv
Note that for any point $x\in\partial\HP\left(v,c_{P,h}(v)\right)\cap\Sph(0,K)$, $x=w+c^2_{P,h}(v)v$ for some $w$ orthogonal to $v$. Moreover, $\|x-c_{P,h}(v)v\|^2=K^2-c^2_{P,h}(v)$ and $\|x-x_n\|^2=(n-c_{P,h}(v))^2+K^2-c^2_{P,h}(v)$ for $x_n=nv$. Thus, we get that:
\[\B\left(x_n,n-c_{P,h}(v)\right)\subset\HP\left(v,c_{P,h}(v),\right)\]
and
\[\HP\left(v,c_{P,h}(v)\right)\cap\supp(P)\subset\B\left(x_n,\sqrt{K^2-c^2_{P,h}(v)+(n-c_{P,h}(v))^2}\right).\]
In particular, since $P\left(\HP\left(v,c_{P,h}(v)\right)\right)<h$, $\B\left(x_n,n-c_{P,h}(v)\right)\subset\B(x_n,\delta_{P,h}(x))$ and since $P\left(\overline{\HP}\left(v,c_{P,h}(v)\right)\right)\geq h$, $\B(x_n,\delta_{P,h}(x))\subset\B\left(x_n,\sqrt{K^2-c^2_{P,h}(v)+(n-c_{P,h}(v))^2}\right)$, with $\delta_{P,h}(x)$ the pseudo-distance defined in Section \ref{section def dtm}.
Finally, for all $y\in\R_d$, if $\langle y,v\rangle=c_{P,h}(v)-\epsilon$ for some $\epsilon>0$, then $\|y-x_n\|^2=\|y\|^2+n^2-2n(c_{P,h}(v)-\epsilon)$, which is superior to $K^2+(n-c_{P,h}(v))^2-c^2_{P,h}(v)$ for $n$ large enough. Thus, for all $n$ large enough, $y\notin\B(x_n,\delta_{P,h}(x_n))$. If $\langle y,v\rangle=c_{P,h}(v)+\epsilon$ for some $\epsilon>0$, then $\|y-x_n\|^2=\|y\|^2+(n-c_{P,h}(v))^2-c^2_{P,h}(v)-2n\epsilon$ which is inferior to $(n-c_{P,h}(v))^2$ for $n$ large enough. Thus, for all $n$ large enough, $y\in\B(x_n,\delta_{P,h}(x_n))$, which concludes the first part of the Lemma.

          Let $(x_n)_{n\geq0}$ be a sequence in $\R_d$ such that $\lim_{n\rightarrow+\infty}\dd_{\overline{\R}_d}(x_n,v_\infty)=0$, that is such that $\lim_{n\rightarrow+\infty}\|x_n\|=+\infty$ and $\lim_{n\rightarrow+\infty}\frac{x_n}{\|x_n\|}=v$. Then,
          \[\|x_n\|-K\leq\delta_{P,h}(x_n)\leq\|x_n\|+K.\]
          Let $y \in \mathbb{R}^d$. Then,
          \begin{align*}
          \|y-x_n\|^2 - \delta_{P,h}(x_n)^2 & = \|y\|^2 - 2 \left\langle x_n,y \right\rangle + O(\|x_n\|) 
                             = \|x_n\|\left ( \frac{\|y\|^2}{\|x_n\|} - 2 \left\langle \frac{x_n}{\|x_n\|},y \right\rangle + O(1) \right ). 
          \end{align*}
          The notation $y_n=O(\|x_n\|)$ means that $\left(\frac{y_n}{\|x_n\|}\right)_{n\in\N}$ is bounded.         
          Thus, up to a subsequence,
          \[\lim_{n\rightarrow+\infty}\frac{\|y-x_n\|^2 - \delta_{P,h}(x_n)^2}{\|x_n\|}=2c-2\langle v,y\rangle,\]
          for some $c\in\R$.
          We deduce that, for all $y\in\R_d\backslash\partial\HP(v,c)$, 
          \[
          \1_{\B(x_n,\delta_{P,h}(x_n))}(y)\rightarrow\1_{\HP(v,c)}(y).
          \]
          In particular, $P(\HP(v,c))\leq h$ and $P(\overline{\HP}(v,c))\geq h$. Therefore, for $P$-almost $y$, $\1_{\HP(v,c)}(y)=\1_{\HP(v,c_{P,h}(v))}(y)$, the result then holds for $c=c_{P,h}(v)$.
          \epv

\subsection{Proof of Lemma \ref{lm enveloppe convexe}}\label{sec:proof_of_lemma_enveloppe_convexe}
\pv

Recall that a $k$-extreme point $x$ of a convex set $S$ is a point $x$ which lies in the interior of a $k$-dimensional convex set within $S$, but not a $k+1$-dimensional convex set within $S$. We will prove that the set of k-extreme points in $\MM_h(P)$ of $\MM_h(P)$ for $k<d$ is equal to $\{m(P_{x,h})\mid x\in\overline{R}_d\backslash\R_d\}$.
In particular, this will yield that the set $\conv(\MM_h(P))$ is equal to $\bigcap_{v\in\Sph(0,1)}H^c(v,\langle m(P_{v_\infty,h}),v\rangle)$.

Let $v\in\Sph(0,1)$,
then by definition the measure $P_{v_\infty,h}$ is supported on $\overline{\HP}(v,c_{P,h}(v))$, satisfies that $hP_{v_\infty,h}$ is a sub-measure of $P$ and the measures $hP_{v_\infty,h}$ and $P$ coincide on $\HP(v,c_{P,h}(v))$. Note that $P(\overline{\HP}(v,c_{P,h}(v)))\geq h$ and $P(\HP(v,c_{P,h}(v)))\leq h$.

We will denote $\overline{C}(P_{v_\infty,h}) = \langle m(P_{v_\infty,h}),v\rangle$, that is, $\overline{C}(P_{v_\infty,h})=P_{v_\infty,h}\langle u,v\rangle$.

Then, for all $x\in\overline{\R}_d$, 
we decompose any measure $P_{x,h}$ as $P_1+P_2$ with $P_1(B) = P_{x,h}(B\cap H(v,c_{P,h}(v)))$ and $P_2(B) = P_{x,h}(B\cap \HP^c(v,c_{P,h}(v)))$.
Note that $P_1$ is also a sub-measure of $P_{v_\infty,h}$. Set $P'_2 = P_{v_\infty,h} - P_1$.
Then, we have
\begin{align*}
P_{x,h}\langle u,v\rangle &= P_1\langle u,v\rangle + P_2\langle u,v\rangle\\
&= P_{v_\infty,h}\langle u,v\rangle - P'_2\langle u,v\rangle + P_2\langle u,v\rangle\\
&= P_{v_\infty,h}\langle u,v\rangle - \langle m(P'_2),v\rangle P'_2(\R_d) + \langle m(P_2),v\rangle P_2(\R_d)\\
&= \bar C(P_{v_\infty,h}) - P'_2\langle u,v\rangle + P_2\langle u,v\rangle\\
&\leq \bar C(P_{v_\infty,h}) - c_{P,h}(v) P'_2(\R_d) + c_{P,h}(v) P_2(\R_d)\\
&= \bar C(P_{v_\infty,h}),
\end{align*}
since $P_2$ is supported on $\HP^c(v,c_{P,h}(v))$ and $P'_2$ is supported on $\overline{\HP}(v,c_{P,h}(v))$ and $P_2(\R_d)=P'_2(\R_d)$.

Thus, for all $x\in\R_d$, $\langle m(P_{x,h}),v\rangle\leq \bar C(P_{v_\infty,h}) = \langle m(P_{v_\infty,h}),v\rangle$. It means that $m(P_{v_\infty,h})$ is not included in a $d$-dimensional simplex within $\conv(\MM_h(P))$. It is thus a k-extreme point of $\conv(\MM_h(P))$ for some $k<d$. Moreover, the hyperplane $\partial H(v,c_{P,h}(v))$ separates $m(P_{v_\infty,h})$ from $\conv(\MM_h(P))$.\\

If $x$ is extreme, then there is some vector $v$ and some constant $C_x$ such that $\langle m(P_{x,h}),v\rangle=C_x$ and such that for all $y$, $\langle m(P_{y,h}),v\rangle\leq C_x$. We aim at proving that $x$ is in $\overline{\R}_d\backslash\R_d$.
Similarly, we get
\begin{align*}
C_x &= P_{x,h}\langle u,v\rangle\\
&= P_1\langle u,v\rangle + P_2\langle u,v\rangle\\
&= P_{v_\infty,h}\langle u,v\rangle - P'_2\langle u,v\rangle + P_2\langle u,v\rangle\\
&\leq C_x - P'_2\langle u,v\rangle + P_2\langle u,v\rangle\\
&\leq C_x - c_{P,h}(v) P'_2(\R_d) + c_{P,h}(v) P_2(\R_d)\\
&= C_x.
\end{align*}

Thus the inequalities are equalities and we get that for $P_2$-almost all $y$, $\langle y,v\rangle=c_{P,h}(v)$ and for $P'_2$-almost all $y$, $\langle y,v\rangle=c_{P,h}(v)$. Thus, $P_{x,h}$ belongs to $\PPP_{v_\infty,h}(P)$. Note that, since there is equality, $C_x=\langle m(P_{v_\infty,h}),v \rangle=\langle m(P_{x,h}),v\rangle$.\\

Note that according to the Krein-Milman theorem, we get that $\conv(\MM_h(P))=\conv(\{m(P_{v_\infty,h})\mid\,v\in\Sph(0,1)\})$.\\

We proved that for all $y\in\R_d$, for all $v\in\Sph(0,1)$, \[\langle m(P_{-v_\infty,h}),v\rangle\leq\langle m(P_{y,h}),v\rangle\leq\langle m(P_{v_\infty,h}),v\rangle.\] Therefore,  the convex set $\conv(\MM_h(P))$ is included in $\bigcap_{v\in\Sph(0,1)} H^c(v,\langle m(P_{v_\infty,h}),v\rangle)$. With the Hahn-Banach separation theorem, we prove that for any $\theta\notin\conv(\MM_h(P))$, there is some vector $v$ such that for all $\theta'\in\conv(\MM_h(P))$, $\langle \theta',v\rangle\leq C<\langle\theta,v\rangle$. In particular, we get that $\langle\theta,v\rangle>\langle m(P_{v_\infty,h}),v\rangle$, meaning that $\theta$ does not belong to $H(v,\langle m(P_{v_\infty,h}),v\rangle)$.
\epv

\subsection{Proof of Lemma \ref{inclusion of conv(M)}}\label{sec:proof_of_lemma_inclusion_of_conv(M)}
\pv
Thanks to Lemma \ref{lm enveloppe convexe}, we have:
\[\conv(\MM_h(P))=\bigcap_{v\in\Sph(0,1)}\HP^c(v,\langle m(P_{v_\infty,h}),v\rangle).\]
Let $0<h'\leq h\leq1$, in order to prove that the map $h\mapsto\conv(\MM_h(P))$ is non-increasing, it is sufficient to prove that
\[\HP^c(v,\langle m(P_{v_\infty,h'}),v\rangle)\supset\HP^c(v,\langle m(P_{v_\infty,h}),v\rangle).\]
Thus, it is sufficient to prove that
\[\langle m(P_{v_\infty,h'}),v\rangle\geq\langle m(P_{v_\infty,h}),v\rangle.\]
Set $P_0$ the sub-measure of $P$ supported on $\overline{\HP}^c(v,c_{P,h}(v))\backslash\HP^c(v,c_{P,h'}(v))$ such that $hP_{v_\infty,h}=h'P_{v_\infty,h'}+(h-h')P_0$.
Then, we have:
\[\langle m(P_{v_\infty,h}),v\rangle=\frac{h'}{h}\langle m(P_{v_\infty,h'}),v\rangle+\frac{h'-h}{h}\langle m(P_0),v\rangle.\]
The results comes from the fact that $\langle m(P_0),v\rangle\leq c_{P,h'}(v)\leq\langle m(P_{v_\infty,h'}),v\rangle$.
\epv
\subsection{Proof of Lemma \ref{lem:transition_M(Pn)_M(P)}}\label{sec:proof_of_lemma_transition_M(Pn)_M}
The proof of Lemma \ref{lem:transition_M(Pn)_M(P)} is based on the following concentration argument, that allows to connect empirical sub-measures with sub-measures for $P_n$. For sake of concision the statement also encompasses sub-Gaussian measures.

\lm\label{lem_concentrationmoy}
           Suppose that $Q\in\PPP^{(V)}(\R_d)$. Then, for every $p>0$, with probability larger than $1-8n^{-p}$, we have,
           \begin{align*}
           \sup_{x,r} | (Q_n - Q)| \1_{\B(x,r)}(y)dy & \leq C \sqrt{\frac{d+1}{n}} + \sqrt{\frac{2p \log(n)}{n}} \\
           \sup_{v,t} | (Q_n - Q)| \1_{\left \langle y , v \right\rangle \leq t}(y)dy & \leq C\sqrt{\frac{d+1}{n}} + \sqrt{\frac{2p\log(n)}{n}} \\
           \sup_{x,r} \| (Q_n - Q) y \1_{\B(x,r)}(y)dy \| & \leq C V \sqrt{d} \frac{(p+1)\log(n)}{\sqrt{n}} \\
           \sup_{v,t} \| (Q_n - Q) y \1_{\left\langle v , y \right\rangle \leq t}dy \| & \leq C V \sqrt{d} \frac{(p+1)\log(n)}{\sqrt{n}} \\
           \sup_{x,r}  \left | (Q_n - Q) \|y\|^2 \1_{\B(x,r)}(y)dy \right | & \leq C V^2 \sqrt{d}\frac{(p+1)\log(n)^{\frac{3}{2}}}{\sqrt{n}} \\
           \sup_{v,t}  \left | (Q_n - Q) \|y\|^2 \1_{\left\langle v , y \right\rangle \leq t}dy \right | & \leq C V^2 \sqrt{d}\frac{(p+1)\log(n)^{\frac{3}{2}}}{\sqrt{n}},          
           \end{align*}  
           where $C>0$ denotes a universal constant.         
\elm
The proof of Lemma \ref{lem_concentrationmoy} is postponed to the following Section \ref{sec:proof_lemma_concentrationmoy}. A significant part of the proof of Lemma \ref{lem:transition_M(Pn)_M(P)} is based on the characterization of $\conv(M_h(P))$ through $\omega^2_{P,h}$ stated by Lemma \ref{lem:omega_convM}, where we recall that $\omega^2_{P,h}(\tau)$ is defined in Definition \ref{def:Opt} by
\begin{align*}
\omega^2_{P,h}(\tau) = \sup_{ x \in \R_d}{\dd^2_{P,h}(x) - \|x-\tau\|^2}.
\end{align*} 
\lm\label{lem_caract_cvx}
Let $C$ denote a convex set, $\theta \in \R_d$, and $\Delta = \dd(\theta,C)$. There exists $v\in\R_d$ with $\|v\|=1$ such that, for all $\tau$ in $C$, 
\[
\left\langle v , \theta - \tau \right\rangle \geq \Delta.
\] 
\elm
\pv[Proof of Lemma \ref{lem_caract_cvx}]
        Denote by $\pi$ the projection onto $C$, and $t=\pi(\theta)$. Then, let $x = \frac{\theta - t}{\Delta}$. We may write
        \begin{align*}
        \left\langle x , \theta -\tau \right\rangle & = \left \langle x , \theta - t \right\rangle + \left\langle x , t - \tau \right\rangle \\
        & = \Delta + \frac{1}{\Delta} \left\langle \theta - t , t-\tau \right\rangle.
        \end{align*}
        Since, for all $\tau$ in $C$, $\left\langle \theta - t , \tau-t \right\rangle \leq 0$, the result follows.
\epv

We are now in position to prove Lemma \ref{lem:transition_M(Pn)_M(P)}.
\pv[Proof of Lemma \ref{lem:transition_M(Pn)_M(P)}]
Let $P$ in $\PPP^{K}(\R_d)$ which puts no mass on hyperplanes nor on spheres, and $\theta\in\conv(\MM_h(P))$. If we choose $p$ large enough (for instance $p=10$), a union bound ensures that the inequalities of Lemma \ref{lem_concentrationmoy} are satisfied for all $n \in \N$ with probability $>0$. Since $P$ puts no mass on hyperplanes, the probability that $n$ points are not in general position is $0$. Hence there exists an empirical distribution $P_n$, in general position, satisfying the inequalities of Lemma \ref{lem_concentrationmoy} for all $n$. In particular, for such a distribution $P_n$  and $(y,r)$ such that $P(\B(y,r)) = h$, we have
           \begin{align*}
          \left \| \frac{P_n u \1_{\B(y,r)}(u)}{P(\B(y,r))} - \frac{P_n u \1_{\B(y,r)}(u)}{P_n(\B(y,r))}\right \| & \leq \frac{K \alpha_n}{h}, \\
        \left \| \frac{P u \1_{\B(y,r)}(u)}{P(\B(y,r))} - \frac{P_n u\1_{\B(y,r)}(u)}{P(\B(y,r))} \right \| & \leq \frac{K \alpha_n}{h}, \\
        \left |(P_n-P) \B(y,r) \right | & \leq \alpha_n,
           \end{align*}
for $\alpha_n \rightarrow 0$. Note that the same holds for means on half-spaces. Now let $x \in \R_d$,                      
\begin{align*}
           \dd^2_{P,h}(x) - \|x-\theta\|^2 &= P_{x,h}\|x-u\|^2 - \|x-\theta\|^2\\
           &=\inf_{y\in\overline{\R}^d}P_{y,h}\|x-u\|^2 - \|x-\theta\|^2\\
           &\geq \inf_{y\in\overline{\R}^d}\|m(P_{y,h})\|^2+v(P_{y,h})-\|\theta\|^2+\inf_{y\in\overline{\R}^d}2\langle x,\theta-m(P_{y,h})\rangle\\
           &\geq-\|\theta\|^2+\inf_{y\in\overline{\R}^d}2\langle x,\theta-m(P_{y,h})\rangle.
           \end{align*}          
           Thus, we may write
           \begin{align*}           
          \inf_{y\in\overline{\R}^d}2\langle x,\theta-m(P_{y,h})\rangle & \geq \min \left [\inf_{y,r|P_n(\B(y,r)) \in \left [ h - \alpha_n,h + \alpha_n \right ]} 2 \left\langle x, \theta - \frac{P_n u\1_{\B_{y,r}(u)}}{P_n(\B(y,r))}  \right\rangle, \right . \\
          & \left .\inf_{v,t|P_n(H(v,t)) \in \left [ h - \alpha_n,h + \alpha_n \right ]} 2 \left\langle x, \theta - \frac{P_n u\1_{H(v,t)(u)}}{P_n(H(v,t))}  \right\rangle \right ] -  \frac{4K \alpha_n \|x\|}{h}, \\
           &=\inf_{\tau\in\bigcup _{s \in \left [ h - \alpha_n,h + \alpha_n \right ]} \MM_s(P_n)} 2 \left\langle x, \theta - \tau \right\rangle - \frac{4K \alpha_n \|x\|}{h}.
           \end{align*}
           Now, if $\dd\left(\theta,\conv\left (\bigcup _{s \in \left [ h - \alpha_n,h + \alpha_n \right ]}  \MM_s(P_n) \right ) \right) = \Delta > \frac{2K \alpha_n}{h}$, then according to Lemma \ref{lem_caract_cvx}, we can choose $x$ in $\R_d$ such that, for all $\tau \in  \conv\left ( \bigcup _{s \in \left [ h - \alpha_n,h + \alpha_n \right ]} \MM_s(P_n) \right )$, 
           \[
           \left\langle \frac{x}{\|x\|}, \theta - \tau \right\rangle - \frac{2K \alpha_n}{h} > 0.
           \]
           In this case, we immediately get
           $\omega^2_{P,h}(\theta) = \sup_{x \in \R_d}\dd^2_{P,h}(x) - \|x-\theta\|^2 = + \infty.$ According to Lemma \ref{lem:omega_convM}, this contradicts $\omega \in \conv(\MM_h(P))$. 

Set $h_n=\frac{q_n}{n}$ for $q_n\in[\![1,n]\!]$ such that $h-\alpha_n\geq h_n\geq h-\alpha_n-\frac{1}{n}$. Note that for $n$ large enough, $h-\alpha_n-\frac{1}{n}>0$, thus $h_n$ is well defined.
Then, according to Lemma \ref{inclusion of conv(M)} and \ref{M(Pn) convex},
\[\conv\left(\bigcup_{s\in[h-\alpha_n,h+\alpha_n]}\MM_s(P_n)\right)\subset\conv\left(\bigcup_{s\in[h_n,1]}\MM_s(P_n)\right)=\MM_{\frac{q_n}{n}}(P_n).\]
Thus, we can build a sequence $(y_n)_{n\geq N}$ for some $N\in\N$ such that $y_n\in\MM_{h_n}(P_n)$ and $\|\theta-y_n\|\leq2\frac{K\alpha_n}{h}$. Hence the result of Lemma \ref{lem:transition_M(Pn)_M(P)}.
\epv

\subsection{Proof of Lemma \ref{lem_concentrationmoy}}\label{sec:proof_lemma_concentrationmoy}
           \pv[Proof of Lemma \ref{lem_concentrationmoy}]
           
           The first inequality is a direct application of Theorem 3.2 in \cite{Boucheron05}, since the Vapnik dimension of balls in $\mathbb{R}^d$ is $d+1$.
           The same argument holds for the second inequality.
           
           Now turn to the third one. Let $\lambda = p \log(n)$, $t =\sqrt{4V^2(\log(n) + \lambda)}$.  Since $Q\in\PPP^{(V)}(\R_d)$, we have that
           \[
           \mathbb{P} \left \{ \max_{i} \|X_i \| \geq t \right \} \leq n e^{-\frac{t^2}{2V^2}} \leq n^{-2p+1}.
           \]
           We may write
           \begin{align*}
           \sup_{x,r} \| (Q_n - Q) y \1_{\B(x,r)}(y)dy \| &= \sup_{x,r} \left \| \frac{1}{n}\sum_{i=1}^{n} {X_i \1_{\B(x,r)}(X_i) - \mathbb{E}(X\1_{\B(x,r)}(X))} \right \| \\
           & \leq \sup_{x,r} \left \| \frac{1}{n}\sum_{i=1}^{n} {X_i \1_{\B(x,r)}(X_i)\1_{\|X_i\| \leq t} - \mathbb{E}(X\1_{\B(x,r)}(X)\1_{\|X\| \leq t})} \right \| \\ 
            & \qquad \qquad + \mathbb{E}(\|X\|\1_{\|X\|>t}) + \sup_{x,r} \frac{1}{n} \sum_{i=1}^{n} \|X_i\|\1_{\|X_i\| >t}.
           \end{align*}
           On one hand, 
           \begin{align*}
           \mathbb{E}(\|X\|\1_{\|X\|>t}) & \leq \sqrt{\mathbb{E}(\|X\|^2)}\sqrt{\mathbb{P}(\|X\| \geq t)} \\
                                        & \leq 2 V e^{-\frac{t^2}{4V^2}} \\
                                        & \leq 2 V n^{-(p+1)}.
           \end{align*}
           On the other hand, with probability larger than $1-n^{-2p+1}$, it holds
           \[
           \sup_{x,r} \frac{1}{n} \sum_{i=1}^{n} \|X_i\|\1_{\|X_i\| >t} =0.
           \]
           Now denote by $f_{x,r,v}$ the function $ \left\langle y \1_{\B(x,r)}(y),v \right\rangle \1_{\|y\| \leq t} $, for $v \in \B(0,1)$, so that 
           \[
           \sup_{x,r} \left \| \frac{1}{n}\sum_{i=1}^{n} {X_i \1_{\B(x,r)}(X_i)\1_{\|X_i\| \leq t} - \mathbb{E}(X\1_{\B(x,r)}(X)\1_{\|X\| \leq t})} \right \| = \sup_{x,r,v} | (Q_n - Q) f_{x,r,v} |.
           \]            
             A straightforward application of MacDiarmid's inequality (see, e.g., \cite[Theorem 6.2]{Massart16}) entails
           \begin{align*}
           \mathbb{P} \left ( \sup_{x,r,v} | (Q_n - Q) f_{x,r,v} | \geq \mathbb{E} \sup_{x,r,v} | (Q_n - Q) f_{x,r,v} | + t \sqrt{\frac{2\lambda}{n}} \right ) \leq e^{-\lambda} = n^{-p}.
           \end{align*}
          It remains to bound  $\mathbb{E} \sup_{x,r,v} | (Q_n - Q) f_{x,r,v} |$. A symmetrization inequality (see, e.g., \cite[Lemma 11.4]{Massart16}) leads to 
          \begin{align*}
          \mathbb{E} \sup_{x,r,v} | (Q_n - Q) f_{x,r,v} | \leq \frac{2t}{n} \mathbb{E}_{X_{1:n}} \mathbb{E}_{\varepsilon} \sup_{x,r,v} \sum_{i=1}^{n}{\varepsilon_i f_{x,r,v}(X_i)/t},
\end{align*}    
where the $\varepsilon_i$'s are i.i.d. Rademacher random variable, and $\mathbb{E}_Y$ denotes expectation with respect to the random variable $Y$.                  
           
           Now suppose that $X_1, \hdots, X_n$ is fixed. In order to apply Dudley's entropy integral we have to provide an upper bound on the metric entropy of $ \mathcal{F} = \left \{ f_{x,r,v}/t \right \}_{x,r,v}$, for the $L_2(P_n)$ distance, that is $d^2(f,f') = \sum_{i=1}^{n} (f(X_i) - f'(X_i))^2/n$. Denote, for any subset of functions $G$, $\mathcal{N}( G,\varepsilon,L_2(P_n))$ the $\varepsilon$-covering number of $G$ with respect to the metric $L_2(P_n)$. Then, if $G \subset G_1 \times G_2$ (for the multiplication), with $\|G_2\| \leq 1$ and $\|G_1\| \leq 1$, we may write
           \begin{align*}
           \mathcal{N}(G,\varepsilon,L_2(P_n)) \leq \mathcal{N}(G_1, \varepsilon/(2\sqrt{2}),L_2(P_n)) \times \mathcal{N}(G_2,\varepsilon/(2\sqrt{2}),L_2(P_n)).  
           \end{align*}
           Define $G_1 = \left \{ (\left\langle v,. \right \rangle)/ t \right \}_{\|v\| \leq 1}$, and $G_2 = \left \{ \1_{\B(x,r) \cap B(0,t)} \right \}_{x,r}$. It is obvious that $\mathcal{F} \subset G_1 \times G_2$. Using Theorem 1 in \cite{Mendelson03}, if $\|G\| \leq 1$, we have
           \[
           \mathcal{N}(G,\varepsilon,L_2(P_n)) \leq \left (\frac{2}{\varepsilon} \right )^{C d_p(G)},
           \]
           where $C$ is an absolute constant and $d_p$ denotes the pseudo-dimension. Hence we have
           \begin{align*}
           \mathcal{N}(G_1,\varepsilon,L_2(P_n)) &\leq \left (\frac{2}{\varepsilon} \right )^{C d} \\
           \mathcal{N}(G_2,\varepsilon,L_2(P_n))&\leq \left (\frac{2}{\varepsilon} \right )^{C (2(d+1))}.
           \end{align*}
We may then deduce
\[
\mathcal{N}(\mathcal{F},\varepsilon,L_2(P_n)) \leq \left ( \frac{4\sqrt{2}}{\varepsilon} \right )^{C(3d+2)}.
\]
Now, using Dudley's entropy integral (see, e.g., \cite[Corollary 13.2]{Massart16}) yields
\begin{align*}
\mathbb{E}_{\varepsilon} \frac{1}{\sqrt{n}}\sup_{x,r,v} \sum_{i=1}^{n}{\varepsilon_i f_{x,r,v}(X_i)/t} & \leq 12 \int_{0}^{1} \sqrt{\log(\mathcal{N}(\mathcal{F},\varepsilon,L_2(P_n)} d\varepsilon \\ 
& \leq 12 \sqrt{C(3d+2)} \int_{0}^{1} \sqrt{\log(\frac{4 \sqrt{2}}{\varepsilon})} d \varepsilon \\
 & \leq C \sqrt{d}. 
\end{align*}                       
            Hence we deduce that
            \[
            \mathbb{E} \sup_{x,r,v} | (Q_n - Q) f_{x,r,v} | \leq \frac{Ct \sqrt{d}}{\sqrt{n}}\leq \frac{CV \sqrt{(p+1)\log(n)}}{\sqrt{n}}.
            \]
            Combining the different terms gives, with probability larger than $1-2n^{-p}$,
            \[
            \sup_{x,r} \| (Q_n - Q) y \1_{\B(x,r)}(y)dy \|  \leq C V \frac{(p+1)\log(n)}{\sqrt{n}}.
            \]
            the third deviation bound follows. The fourth deviation bound may be proved the same way. 
            
            For the $5$-th inequality, as before let $\lambda = p \log(n)$ and $t=\sqrt{4V^2(\log(n) + \lambda)}$. Similarly, we may write
             \begin{align*}
           \sup_{x,r} | (Q_n - Q) \|y\|^2 \1_{\B(x,r)}(y)dy \| &= \sup_{x,r} \left | \frac{1}{n}\sum_{i=1}^{n} {\|X_i\|^2 \1_{\B(x,r)}(X_i) - \mathbb{E}(\|X\|^2\1_{\B(x,r)}(X))} \right \| \\
           & \leq \sup_{x,r} \left | \frac{1}{n}\sum_{i=1}^{n} {\|X_i\|^2 \1_{\B(x,r)}(X_i)\1_{\|X_i\| \leq t} - \mathbb{E}(\|X\|^2\1_{\B(x,r)}(X)\1_{\|X\| \leq t})} \right | \\ 
            & \qquad \qquad + \mathbb{E}(\|X\|^2\1_{\|X\|>t}) + \sup_{x,r} \frac{1}{n} \sum_{i=1}^{n} \|X_i\|^2\1_{\|X_i\| >t},
           \end{align*}
           with $\mathbb{E}(\|X\|^2\1_{\|X\|>t}) \leq 2V^2n^{-(p+1)}$ and, with probability larger than $1-n^{-\frac{1}{p}}$, 
            \[
            \sup_{x,r} \frac{1}{n} \sum_{i=1}^{n} \|X_i\|^2\1_{\|X_i\| >t} =0.
            \]
            Using \cite[Theorem 6.2]{Massart16} again leads to 
            \begin{multline*}
            \mathbb{P} \left (\sup_{x,r} \left | \frac{1}{n}\sum_{i=1}^{n} {\|X_i\|^2 \1_{\B(x,r)}(X_i)\1_{\|X_i\| \leq t} - \mathbb{E}(\|X\|^2\1_{\B(x,r)}(X)\1_{\|X\| \leq t})} \right | \right . \\ 
            \geq \left . \mathbb{E} \sup_{x,r} \left | \frac{1}{n}\sum_{i=1}^{n} {\|X_i\|^2 \1_{\B(x,r)}(X_i)\1_{\|X_i\| \leq t} - \mathbb{E}(\|X\|^2\1_{\B(x,r)}(X)\1_{\|X\| \leq t})} \right | + t^2\sqrt{\frac{2\lambda}{n}} \right ) \leq e^{-\lambda} =n^{-\frac{1}{p}}.
            \end{multline*}
            At last, combining a symmetrization inequality with a contraction principle (\cite[Theorem 11.5]{Massart16}) gives
     \begin{align*}
     \mathbb{E} & \sup_{x,r} \left | \frac{1}{n}\sum_{i=1}^{n} {\|X_i\|^2 \1_{\B(x,r)}(X_i)\1_{\|X_i\| \leq t} - \mathbb{E}(\|X\|^2\1_{\B(x,r)}(X)\1_{\|X\| \leq t})} \right | \\
     & \leq \frac{2t^2}{n} \mathbb{E}_{X_{1:n}} \mathbb{E}_{\varepsilon} \sup_{x,r} \sum_{i=1}^{n} \varepsilon_i \frac{\|X_i\|^2}{t^2} \1_{\B(x,r)\cap \B(0,t)}(X_i) \\
     & \leq \frac{2t^2}{n} \mathbb{E}_{X_{1:n}} \mathbb{E}_{\varepsilon} \sup_{x,r} \sum_{i=1}^{n} \varepsilon_i \1_{\B(x,r)\cap \B(0,t)}(X_i) \\
     & \leq \frac{C\sqrt{d}t^2}{\sqrt{n}},
\end{align*}            
where the last line may be derived the same way as for the third inequality, combining \cite[Theorem 1]{Mendelson03} and \cite[Corollary 13.2]{Massart16}. Gluing all pieces yields, with probability larger than $1-2n^{-\frac{1}{p}}$,
\[
\sup_{x,r}  \left | (Q_n - Q) \|y\|^2 \1_{\B(x,r)}(y)dy \right |  \leq C V^2 \sqrt{d}\frac{(p+1)\log(n)^{\frac{3}{2}}}{\sqrt{n}}.
\]
The last inequality follows from the same argument.
           \epv


\subsection{Proof of Lemma \ref{lemme de concavite}}
\pv
In Proposition 3.6 from \cite{Merigot1}, we get:
\begin{align*}
\dd^2_{P,h}(y)&=P_{y,h}\|y-u\|^2\\
&\leq P_{x,h}\|y-u\|^2\\
&=\|y-x\|^2+P_{x,h}\|x-u\|^2+2\langle y-x,x-P_{x,h}u\rangle.
\end{align*}
In particular,
\[\dd^2_{P,h}(y)-\|y\|^2\leq\dd^2_{P,h}(x)-\|x\|^2-2\langle y-x,P_{x,h}u\rangle,\]
with equality if and only if $P_{y,h}\|y-u\|^2=P_{x,h}\|y-u\|^2$, that is if and only if like $hP_{y,h}$, $hP_{x,h}$ is also a sub-measure of $P$ with total mass $h$, whose support is contained in the closed ball $\overline{\B}(y,\delta_{P,h}(y))$ and whose restriction to the open ball $\B(y,\delta_{P,h}(y))$ coincides with $P$; see \cite{Merigot1}, Proposition 3.3.
\epv

\section{Proofs for Section \ref{sec:kPDTM_a_coreset_for_DTM}}
\subsection{Proof of Theorem \ref{opt nonempty}}\label{sec:proof_thm_opt_nonempty}
\pv
First note that for all $t,s\in\overline{\R}_d^{(k)}$, denoting by
\[f_s:x\mapsto\min_{i\in[\![1,k]\!]}\left(\|x-m(P_{s_i,h})\|^2+v(P_{s_i,h})\right),\] we have:
\[Pf_s(u)-f_t(u)\leq\sum_{i=1}^k\tilde P_{t_i,h}(\R_d)\left(2\langle\tilde m(\tilde P_{t_i,h}),m(P_{t_i,h})-m(P_{s_i,h})\rangle+M(P_{s_i,h})-M(P_{t_i,h})\right).\]
Then, according to Lemma \ref{lem_cv_ball_plane} and the dominated convergence Theorem, for any $v\in\Sph(0,1)$, there is a sequence $(x_n)_{n\in\N}$ in $\R_d$ such that $m(P_{x_n,h})\rightarrow m(P_{v_\infty,h})$ and $M(P_{x_n,h})\rightarrow M(P_{v_\infty,h})$. Then, $\lim\sup_{n\rightarrow+\infty}Pf_{x_n}(u)-f_{v_\infty}(u)\leq0$. Thus, $\inf_{t\in\R_d}Pf_t(u)=\inf_{t\in\overline{\R}_d}Pf_t(u)$.

Let $(t_n)_{n\in\N}$ be a sequence in $\R_d^{(k)}$ such that $Pf_{t_n}(u)\leq\inf_{t\in\overline{\R}_d}Pf_t(u)+\frac{1}{n}$, and denote by $m^*$ the limit of a converging subsequence of $\left(\tilde m(\tilde P_{t_{n,1},h}),\tilde m(\tilde P_{t_{n,2},h}),\ldots,\tilde m(\tilde P_{t_{n,k},h})\right)_{n\in\N}$ in the compact space $\overline{\B}(0,K)^{(k)}$.
Then, thanks to Lemma \ref{lemme de concavite}, and recalling that $\forall y\in\R_d,\,\dd^2_{P,h}(y)=\|y-m(P_{y,h})\|^2+v(P_{y,h})$,
\[Pf_{m^*}(u)-f_{t_n}(u)\leq2\sum_{i=1}^k\tilde P_{t_{n,i},h}(\R_d)\langle \tilde m(\tilde P_{t_{n,i},h})-m^*_i,m(P_{t_{n,i},h})-m(P_{m^*_i,h})\rangle,\]
which goes to zero when $n\rightarrow+\infty$ since $\|m(P_{y,h})\|\leq K$ whenever $y \in \R_d$. Thus, $Pf_{m^*}(u)=\inf_{t\in\overline{\R}_d}Pf_t(u)$.
In particular, there is some $s\in\opt(P,h,k)\cap\overline{\B}(0,K)^{(k)}$. Take $P_{s_i,h}\in\PPP_{s_i,h}(P)$ for all $i\in[\![1,k]\!]$, then set $s^*\in\overline{\B}(0,K)^{(k)}$ such that $s^*_i=\tilde m(\tilde P_{s_{i},h})$ for all $i\in[\![1,k]\!]$. Then for any choice of $P_{s^*_i,h}\in\PPP_{s^*_i,h}(P)$,
\begin{align*}
0&\leq Pf_{s^*}(u)-f_{s}(u)\\
&\leq\sum_{i=1}^k\tilde P_{s_{i},h} \|u-m(P_{s^*_i,h})\|^2+v(P_{s^*_i,h})-\|u-m(P_{s_i,h})\|^2-v(P_{s_i,h})\\
&=\sum_{i=1}^k\tilde P_{s_{i},h}\left(\dd^2_{P,h}(s^*_i)-\|s^*_i\|^2\right)-\left(\dd^2_{P,h}(s_i)-\|s_i\|^2\right)-\langle u-s^*_i,2m(P_{s^*_{i},h})\rangle+\langle u-s_i,2m(P_{s_{i},h})\rangle\\
&\leq2\sum_{i=1}^k\tilde P_{s_{i},h}(\R_d)\langle \tilde m(\tilde P_{s_{i},h})-s^*_i,m(P_{s_{i},h})-m(P_{s^*_i,h})\rangle=0.
\end{align*}
Thus, inequalities are all equalities. In particular, equality in Lemma \ref{lemme de concavite} leads to $P_{s_i,h}\in\PPP_{s^*_i,h}(P)$, and by choosing $P_{s^*_i,h}=P_{s_i,h}$, the Laguerre measures $(\tilde P_{s_i,h})_{i\in[\![1,k]\!]}$ are also appropriate for $s^*$. Then, $\tilde m(\tilde P_{s^*_i,h})=\tilde m(\tilde P_{s_i,h})=s^*_i$. Thus, $s^*\in\opt(P,h,k)\cap\overline{\B}(0,K)^{(k)}$ and satisfies for some $(P_{s^*_i,h})_{i\in[\![1,k]\!]}$, $\tilde m(\tilde P_{s^*_i})=s^*_i$, for all $i\in[\![1,k]\!]$.
\epv

\subsection{Proof of Corollary \ref{cor:approx_lebesgue}}\label{sec:proof_cor_approx_lebesgue}
\pv[Proof of Corollary \ref{cor:approx_lebesgue}]
The proof of Corollary \ref{cor:approx_lebesgue} is based on the following bounds, in the case where
 $P$ is absolutely continuous with respect to the Lebesgue measure, with density $f$ satisfying $0<f_{min} \leq f \leq f_{max}$. 
\begin{align}
f_{M}^{-1}(k) & \leq  2K \sqrt{d} k^{-1/d} \label{eq:covering_lebesgue} \\ 
\zeta_{P,h}(f_{M}^{-1}(k)) & \leq K C_{f_{max},K,d,h} k^{-1/d}. \label{eq:meanlipschitz_lebesgue}
\end{align}
      The first equation proceeds from the following. Since $M \subset \B(0,K)$, for any $\varepsilon >0$ we have
      \begin{align*}
      f_M(\varepsilon) \leq f_{\B}(0,K)(\varepsilon) \leq \left ( \frac{2K \sqrt{d}}{\varepsilon} \right )^{d}.
\end{align*}       Hence \eqref{eq:covering_lebesgue}. To prove the second inequality, we will use the following Lemma.
      \lm\label{lem:mean_modulus_continuity}
      Suppose that $P$ has a density $f$ satisfying $0 < f_{min} \leq f \leq f_{max}$. Let $x$, $y$ be in $M$, and denote by $\delta = \|x-y\|$. Then
      \[
      \| m(P_{x,h}) - m(P_{y,h}) \| \leq {(2K)^{d+1} \omega_d} \left ( 1 + \delta \left ( \frac{f_{max} \omega_d}{h} \right )^{1/d} \right )^{d-1} \left( \frac{f_{max} \omega_d}{h} \right )^{1/d} \delta.
      \]      
      \elm
      \pv[Proof of Lemma \ref{lem:mean_modulus_continuity}]     
      Since $P$ has a density,  $P \partial \B(x,\delta_{x,h})=P \partial \B(y,\delta_{y,h})=0$. We deduce that $P_{x,h} = \frac{1}{h}P_{|\B(x,h)}$ and  $P_{y,h} = \frac{1}{h}P_{|\B(y,h)}$. Without loss of generality, assume that $\delta_{x,h} \geq \delta_{y,h}$. Then $\B(y,\delta_{y,h}) \subset \B(x,\delta_{x,h} + \delta)$. We may write
      \begin{align*}
      \| m(P_{x,h}) - m(P_{y,h}) \| & = \frac{1}{m_0} \left \| P \left ( u (\mathbbm{1}_{\B(x,\delta_{x,h})}(u) - \mathbbm{1}_{\B(y,\delta_{y,h})}(u) \right ) \right \| \\
       & \leq \frac{2K}{h}  P \left | (\mathbbm{1}_{\B(x,\delta_{x,h})}(u) - \mathbbm{1}_{\B(y,\delta_{y,h})}(u) \right | \\
       & \leq \frac{2K}{h} P \left ( \B(x, \delta_{x,h} + \delta) \cap \B(x, \delta_{x,h} + \delta)^c \right ) \\
       & \leq \frac{2K}{h} \omega_{d} \left [ (\delta_{x,h} + \delta)^d - \delta_{x,h}^d \right ] \\
       & \leq \frac{(2K)^{d+1} \omega_d}{h} \left [ (1 + \frac{\delta}{\delta_{x,h}})^d -1 \right ].
      \end{align*}
      Since $(1+ v)^d \leq 1+d(1+v)^{d-1}v$, for $v \geq 0$, and $ \delta_{x,h} \geq \left ( \frac{h}{f_{max} \omega_d} \right )^{1/d}$, the result follows.
      
      \epv
Hence \eqref{eq:meanlipschitz_lebesgue}. The result of Corollary \ref{cor:approx_lebesgue} follows.      
\epv
\subsection{Proof of Corollary \ref{cor:approx_manifold}}\label{sec:proof_cor_approx_manifold}
\pv[Proof of Corollary \ref{cor:approx_manifold}]
     Without loss of generality we assume that $N$ is connected. Since $P$ has a density with respect to the volume measure on $N$, we have $P(N^{\mathrm{o}})=1$. Thus we take $M=N^{\mathrm{o}}$, that is the set of interior points. Since $P$ satisfies a $(cf_{min},d')$-standard assumption, we have 
     \[
     f_M(\varepsilon) \leq \frac{2^{d'}}{c f_{min}} r^{-d'},
     \]
     according to \cite[Lemma 10]{Chazal15}. Hence $f_{M}^{-1}(k) \leq C_{f_{min},N} k^{-1/d'}$. It remains to bound the continuity modulus of $x \mapsto m(P_{x,h})$. For any $x$ in $M$, since $P(\partial N) =0$ and $P$ has a density with respect to the volume measure on $N$, we have $P_{x,h} =P_{|\B(x,h)}$. Besides, since for all $r >0$ $P(\B(x,r)) \geq c f_{min} r^{d'}$, we may write $\delta_{x,h} \leq c_{N,f_{min}} h^{1/d'} \leq \rho/12$, for $h$ small enough. Now let $x$ and $y$ be in $M$ so that $\|x-y\| = \delta \leq \rho/12$, and without loss of generality assume that $\delta_{x,h} \geq \delta_{y,h}$. Then, proceeding as in the proof of Lemma \ref{lem:mean_modulus_continuity}, it comes
     \[
     \left \| m(P_{x,h}) - m(P_{y,h}) \right \| \leq \frac{2K}{h} P \left ( \B(x,\delta_{x,h} + \delta) \cap \B(x,\delta_{x,h})^c \right ).
     \]
     Since $\delta_{x,h} + \delta \leq \rho/ 6$, for any $u$ in  $\B(x,\delta_{x,h} + \delta) \cap M$ we may write $u=\exp_{x}(rv)$, where $v \in T_xM$ with $\|v\|=1$ and $r=d_N(u,x)$ is the geodesic distance between $u$ and $x$ (see, e.g., \cite[Proposition 25]{Aamari15}). Note that, according to \cite[Proposition 26]{Aamari15}, for any $u_1$ and $u_2$ such that $\|u_1 - u_2\| \leq \rho/4$, 
     \begin{align}
     \|u_1 - u_2\| \leq d_N(u_1,u_2) \leq 2 \|u_1-u_2\| \label{eq:euclidean_geodesic_distance}. 
     \end{align}
     Now let $p_1, \hdots, p_m$ be a $\delta$-covering set of the sphere $\mathcal{S}_{x,\delta_{x,h}}=\left \{ u \in M | \|x-u\| = \delta_{x,u} \right \}$. According to \eqref{eq:euclidean_geodesic_distance}, we may choose $m \leq c_{d'} \delta_{x,h}^{d'-1} \delta^{-(d'-1)}$.

     Now, for any $u$ such that $u \in M$ and $\delta_{x,h} \leq \|x-u\| \leq \delta_{x,h} + \delta$, there exists $t \in \mathcal{S}_{x,\delta_{x,h}}$ such that $\|t-u\| \leq 2 \delta$. Hence  
     \begin{align*}
     P \left ( \B(x,\delta_{x,h} + \delta) \cap \B(x,\delta_{x,h})^c \right ) & \leq \sum_{j=1}^{m} P \left ( \B(p_j,2\delta) \right ).
     \end{align*}
     
     Now, for any $j$, since $2 \delta \leq \rho/6$, in local polar coordinates around $p_j$ we may write, using \eqref{eq:euclidean_geodesic_distance} again,
     \begin{align*}
     P \left ( \B(p_j,2\delta) \right ) \leq  \int_{r,v|\exp_{p_j}(rv) \in M,r \leq 4 \delta} f(r,v) J(r,v) dr dv \\
     & \leq f_{max} \int_{r,v|r \leq 4 \delta} J(r,v) dr dv
     \end{align*}
where $J(r,v)$ denotes the Jacobian of the volume form. According to \cite[Proposition 27]{Aamari15}, we have $J(r,v) \leq C_{d'} r^{d'}$. Hence $P \left ( \B(p_j,2\delta) \right ) \leq C_{d'}f_{max} \delta^{d'}$. We may conclude
\begin{align*}
\left \| m(P_{x,h}) - m(P_{y,h}) \right \| &\leq \frac{2K}{h} m C_{d'}f_{max} \delta^{d'} \\
& \leq C_{N,f_{max},f_{min}} \delta.
\end{align*}      
Choosing $k$ large enough so that $f_{M}^{-1}(k) \leq C_{f_{min},N} k^{-1/d'} \leq \rho/12$ gives the result of Corollary \ref{cor:approx_manifold}.
\epv
\subsection{Proof of Proposition \ref{prop approximation Wasserstein}}\label{sec:proof_prop_approximation_wasserstein}
\pv
To lighten the notation we omit the $\varepsilon$ in $\dd_{Q,h,k,\varepsilon}$. For all $x\in\supp(P)$,
\begin{align*}
\dd^2_{Q,h,k}(x)-\dd^2_{P,h}(x)&=\dd^2_{Q,h,k}(x)-\dd^2_{Q,h}(x) + \dd^2_{Q,h}(x)-\dd^2_{P,h}(x)\\
&\geq-\|\dd^2_{P,h}-\dd^2_{Q,h}\|_{\infty,\supp(P)}.
\end{align*}
Thus, $\left(\dd^2_{Q,h,k}-\dd^2_{P,h}\right)_-\leq\|\dd^2_{P,h}-\dd^2_{Q,h}\|_{\infty,\supp(P)}$ on $\supp(P)$, where $f_-:x\mapsto f(x)\1_{f(x)\leq0}$ denotes the negative part of any function $f$ on $\R_d$.
Then,
\begin{align*}
P\left|\dd^2_{Q,h,k}-\dd^2_{P,h}\right|(u)&=P \dd^2_{Q,h,k}(u)-\dd^2_{P,h}(u) + 2\left(\dd^2_{Q,h,k}(u)-\dd^2_{P,h}(u)\right)_-\\
&\leq P\Delta(u) + P\dd^2_{P,h,k}(u) - \dd^2_{P,h}(u) + 2\|\dd^2_{P,h}-\dd^2_{Q,h}\|_{\infty,\supp(P)}.
\end{align*}
with $\Delta = \dd^2_{Q,h,k}-\dd^2_{P,h,k}$.
We can bound $P\Delta(u)$ from above. Let $s\in\opt(P,h,k)\cap\overline{\B}(0,K)^{(k)}$ such that $s_i=\tilde m(\tilde P_{s_i,h})$ for all $i\in[\![1,k]\!]$. Such an $s$ exists according to Theorem \ref{opt nonempty}.
Set $f_{Q,t}(x)=2\langle x,m(Q_{t,h})\rangle+v(Q_{t,h})$ for $t\in\R_d$, and let $t\in\opt(Q,h,k)$.
\begin{align*}
P\Delta(u)&=P\min_{i\in[\![1,k]\!]}f_{Q,t_i}(u)-\min_{i\in[\![1,k]\!]}f_{P,s_i}(u)\\
&\leq (P-Q)\min_{i\in[\![1,k]\!]}f_{Q,t_i}(u) + \epsilon + (Q-P)\min_{i\in[\![1,k]\!]}f_{P,s_i}(u) + P \min_{i\in[\![1,k]\!]}f_{Q,t_i}(u)-\min_{i\in[\![1,k]\!]}f_{P,s_i}(u).
\end{align*}
For any transport plan $\pi$ between $P$ and $Q$,
\begin{align*}
P-Q \min_{i\in[\![1,k]\!]}f_{Q,t_i}(u)&=\E_{(X,Y)\sim\pi}\left[\min_{i\in[\![1,k]\!]}2\langle X,m(Q_{t,h})\rangle+v(Q_{t,h})-\min_{i\in[\![1,k]\!]}2\langle Y,m(Q_{t,h})\rangle+v(Q_{t,h})\right]\\
&\leq2\E_{(X,Y)\sim\pi}\left[\sup_{t\in\overline{\R}_d}\langle X-Y,m(Q_{t,h})\rangle\right].
\end{align*}
Thus, $P-Q \min_{i\in[\![1,k]\!]}f_{Q,t_i}(u)\leq 2W_1(P,Q)\sup_{t\in\overline{\R}_d}m(Q_{t,h})$, after taking for $\pi$ the optimal transport plan for the $L_1$-Wasserstein distance (noted $W_1$) between $P$ and $Q$.

Also note that $P \min_{i\in[\![1,k]\!]}f_{Q,t_i}(u)-\min_{i\in[\![1,k]\!]}f_{P,s_i}(u)$ is bounded from above by
\begin{align*}
&\leq\sum_{i=1}^k\tilde P_{s_i}\left(2\langle u,m(Q_{s_i,h})\rangle+v(Q_{s_i,h})\right)-\left(2\langle u,m(P_{s_i,h})\rangle+v(P_{s_i,h}))\right)\\
&=\sum_{i=1}^k\tilde P_{s_i}\min_{j\in[\![1,k]\!]}\left(2\langle u-s_i,m(P_{s_i,h})-m(Q_{s_i,h})\rangle+\dd^2_{Q,h}(s_i)-\dd^2_{P,h}(s_i)\right)\\
&\leq\|\dd^2_{P,h}-\dd^2_{Q,h}\|_{\infty,\supp(P)}+2\sum_{i=1}^k\tilde P_{s_i,h}(\R_d)\langle\tilde m(\tilde P_{s_i,h})-s_i,m(P_{s_i,h})-m(Q_{s_i,h})\rangle.
\end{align*}
Since $s_i=\tilde m(\tilde P_{s_i,h})$, the result follows.

\epv
\subsection{Proof of Proposition \ref{prop:bound_ktdm_dtosupport}}\label{sec:proof_bound_fdtm_dtosupport}
\pv
The proof of Proposition \ref{prop:bound_ktdm_dtosupport} relies on \cite[Corollary 4.8]{Merigot1}. Namely, if $P$ satisfies \eqref{eq:abstandard}, then
\[
\|d_{P,h} - d_M \|_{\infty} \leq C(P)h^{-\frac{1}{d'}}.
\]
Let $\Delta_{\infty,K}$ denote $\sup_{x \in M}|d_{Q,h,k,\varepsilon}|$, and let $x \in M$ achieving the maximum distance. Since $d_{Q,h,k,\varepsilon}$ is $1$-Lipschitz, we deduce that $\B(x,\Delta_{\infty}/2) \subset \{y| \quad |d_{Q,k,h,\varepsilon}(y)| \geq \Delta_{\infty}/2 \}$.
            Since $P(\B(x,\Delta_{\infty}/2) \geq C(P)\Delta_{\infty}^{d'}$, Markov inequality yields that
      \[
      \Delta_P^2 \geq  C(P)\Delta_{\infty}^{d'+2}.
      \]
Thus we have $\sup_{x \in M} |d_{Q,h,k} - d_M|(x) \leq  C(P)^{-\frac{1}{d'+2}} \Delta_{P}^{\frac{2}{d'+2}}$. Now, for $x \in \mathbb{R}^d$, we let $p \in M$ such that $\|x-p\| =d_M(x)$. Denote by $r=\|x-p\|$, and let $t_j$ be such that $d_{Q,h,k,\varepsilon}(p) = \sqrt{\|p-m(Q_{t_j,h})\|^2 + v(Q_{t_j,h})}$. Then
      \begin{align*}
      d_{Q,h,k}(x) &\leq \sqrt{\|x-m(Q_{t_j,h})\|^2 + v(Q_{t_j,h})} \\
             & \leq \sqrt{d_{Q,h,k}^2(p) + r^2 +2r \|p -m(Q_{t_j,h}) \|} \\
             & \leq \sqrt{d_{Q,h,k}^2(p) + r^2 +2r d_{Q,h,k}(p)} \\
             & \leq r + (d_{Q,h,k}(p) - d_M(p)). 
      \end{align*}    
On the other hand, we have $d_{Q,h,k,\varepsilon} \geq d_{Q,h}$, along with $\|d_{Q,h} - d_{P,h}\|_{\infty} \leq h^{-\frac{1}{2}} W_2(P,Q)$ (see, e.g., \cite[Theorem 3.5]{Merigot1}) as well as $d_{P,h} \geq d_M$. Hence
\[
d_{Q,h,k,\varepsilon} \geq d_M - h^{-\frac{1}{2}} W_2(P,Q).
\]      
\epv
\section{Proofs for Section \ref{sec:Approximation_kdtm_pointclouds}}
\subsection{Proof of Proposition \ref{prop algo}}\label{sec:proof_prop_algo}
\pv
For any $t=(t_1,t_2,\ldots t_k)\in\R_d^{(k)}$, we note $c_i=\frac{\sum_{X\in\CC(t_i)}X}{|\CC(t_i)|}$. Then,
\begin{align*}
&P_n\min_{i\in[\![1,k]\!]}\|u-m(P_{n\,t_i,h})\|^2+v(P_{n\,t_i,h})\\
&=\sum_{i=1}^k\frac{1}{n}\sum_{X\in\CC(t_i)}\|X-m(P_{n\,t_i,h})\|^2+v(P_{n\,t_i,h})\\
&=\sum_{i=1}^k\frac{1}{n}\sum_{X\in\CC(t_i)}\|X\|^2-2\langle X-t_i,m(P_{n\,t_i,h})\rangle+\left(\dd^2_{P_n,h}(t_i)-\|t_i\|^2\right)\\
&=\frac{1}{n}\sum_{j=1}^n\|X_j\|^2+\sum_{i=1}^k\frac{|\CC(t_i)|}{n}\left(-2\langle c_i-t_i,m(P_{n\,t_i,h})\rangle+\left(\dd^2_{P_n,h}(t_i)-\|t_i\|^2\right)\right)\\
&\geq\frac{1}{n}\sum_{j=1}^n\|X_j\|^2+\sum_{i=1}^k\dd^2_{P_n,h}(c_i)-\|c_i\|^2\\
&=\sum_{i=1}^k\frac{1}{n}\sum_{X\in\CC(t_i)}\|X-m(P_{n\,c_i,h})\|^2+v(P_{n\,c_i,h})\\
&\geq\sum_{i=1}^k\frac{1}{n}\sum_{X\in\CC(c_i)}\|X-m(P_{n\,c_i,h})\|^2+v(P_{n\,c_i,h})\\
&=P_n\min_{i\in[\![1,k]\!]}\|u-m(P_{n\,c_i,h})\|^2+v(P_{n\,c_i,h}).
\end{align*}
We used Lemma \ref{lemme de concavite}.
\epv

\subsection{Proof of Theorem \ref{thm:prox_empirical_version}}\label{sec:proof_thm_prox_empirical_version}

Let $\gamma$ and $\hat\gamma$ the functions defined for $(t,x)\in\R_d^{(k)}\times\R_d$ with $t=(t_1,t_2,\ldots,t_k)$, by:
\[\gamma(t,x)=\min_{i\in[\![1,k]\!]}-2\langle x,m(Q_{t_i,h})\rangle+\|m(Q_{t_i,h})\|^2+v(Q_{t_i,h}),\]
and
\[\hat\gamma(t,x)=\min_{i\in[\![1,k]\!]}-2\langle x,m(Q_{n\,t_i,h})\rangle+\|m(Q_{n\,t_i,h})\|^2+v(Q_{n\,t_i,h}).\]
The proof of Theorem \ref{thm:prox_empirical_version} is based on the two following deviation Lemmas.
\lm
\label{lem:Borne Q-Qn}
If $Q$ is sub-Gaussian with variance $V^2$, then, for every $p>0$, with probability larger than $1-2n^{-\frac{1}{p}}$, we have
\[\sup_{t\in\R_d^{(k)}} \left | (Q-Q_n)\gamma(t,u) \right |\leq C \frac{\sqrt{kd} V^2 \log(n)}{h \sqrt{n}}.\]
\elm
The proof of Lemma \ref{lem:Borne Q-Qn} is deferred to Section \ref{sec:proof_lemma_Q_Qn}.
\lm
\label{Borne gamma-hat gamma}
Assume that $Q$ is sub-Gaussian with variance $V^2$, then, for every $p>0$, with probability larger than $1-7n^{-p}$, we have
\[
\sup_{t \in \mathbb{R}_d^{(k)}} \left | Q_n(\gamma-\hat{\gamma})(t,u) \right | \leq C V^2 \frac{(p+1)^{\frac{3}{2}}\log(n)^{\frac{3}{2}}}{h \sqrt{n}}.
\]
\elm
As well, the proof of Lemma \ref{Borne gamma-hat gamma} is deferred to Section \ref{sec:proof_borne_gammahat_gamma}. We are now in position to prove Theorem \ref{thm:prox_empirical_version}.
\pv[Proof of Theorem \ref{thm:prox_empirical_version}]
Let \[s = \argmin\left\{Q\gamma(t,u)\mid t=(t_1,t_2,\ldots t_k)\in\R_d^{(k)}\right\},\]
\[\hat s = \argmin\left\{Q_n\hat\gamma(t,u)\mid t=(t_1,t_2,\ldots t_k)\in\R_d^{(k)}\right\}\]
and
\[\tilde s = \argmin\left\{Q_n\gamma(t,u)\mid t=(t_1,t_2,\ldots t_k)\in\R_d^{(k)}\right\}.\]

With these notations, for all $x\in\R_d$, $\dd^2_{Q,h,k}(x)=\|x\|^2+\gamma(s,x)$ and $\dd^2_{Q_n,h,k}(x)=\|x\|^2+\hat\gamma(\hat s,x)$.
We intend to bound $l(s,\hat s)=Q(\dd^2_{Q_n,h,k}(u)-\dd^2_{Q,h,k}(u))$, which is also equal to $l(s,\hat s)=Q(\gamma(\hat s,u)-Q\gamma(s,u))$.

We have that:
\begin{align*}
l(s,\hat s)&=Q\gamma(\hat s,u)-Q_n\gamma(\hat s,u)+Q_n\gamma(\hat s,u)-Q_n\gamma(\tilde s,u)+Q_n\gamma(\tilde s,u)-Q\gamma(s,u)\\
&\leq \sup_{t\in\R_d^{(k)}}(Q-Q_n)\gamma(t,u)+Q_n(\gamma-\hat\gamma)(\hat s,u) \\
&+Q_n(\hat\gamma(\hat s,u)-\hat\gamma(\tilde s,u))+Q_n(\hat\gamma-\gamma)(\tilde s,u)+\sup_{t\in\R_d^{(k)}}(Q_n-Q)\gamma(t,u),
\end{align*}
where we used the fact that $Q_n\gamma(\tilde s,u)\leq Q_n\gamma(s,u)$.
 Now, since $Q_n(\hat\gamma(\hat s,u)-\hat\gamma(\tilde s,u))\leq 0$, we get:
\begin{align*}
l(s,\hat s)&\leq\sup_{t\in\R_d^{(k)}}(Q-Q_n)\gamma(t,u)+\sup_{t\in\R_d^{(k)}}(Q_n-Q)\gamma(t,u)\\
&+\sup_{t\in\R_d^{(k)}}Q_n(\gamma-\hat\gamma)(t,u)+\sup_{t\in\R_d^{(k)}}Q_n(\hat\gamma-\gamma)(t,u).
\end{align*}
Combining Lemma \ref{lem:Borne Q-Qn} and Lemma \ref{Borne gamma-hat gamma} entails, with probability larger than $1-8n^{-p}$, 
\begin{align*}
l(s,\hat s)&\leq C \left (V^2 \frac{(p+1)^{\frac{3}{2}}\log(n)^{\frac{3}{2}}}{h \sqrt{n}} +\frac{\sqrt{kd} V^2 \log(n)}{h \sqrt{n}} \right ).
\end{align*}
It remains to bound $|P\dd^2_{Q_n,h,k} - Q\dd^2_{Q_n,h,k}|$ as well as $|P\dd^2_{Q,h,k} - Q\dd^2_{Q,h,k}|$. To this 	aim we recall that $X = Y +Z$, $Z$ being sub-Gaussian with variance $\sigma^2$. Thus, denoting by $t_j(x) = \arg\min_j \|x-m(Q_{t_j,h})\|^2 + v(Q_{t_j,h})$,
\begin{align*}
P\dd^2_{Q,h,k} - Q\dd^2_{Q,h,k} & \leq \mathbb{E} \left [ \|Y-m(Q_{t_j(Y),h})\|^2 + v(Q_{t_j(Y),h}) - \left (\|Y+Z-m(Q_{t_j(Y),h})\|^2 + v(Q_{t_j(Y),h}) \right ) \right ]  \\
                   & \leq \mathbb{E} \|Z\|^2 + 2 \mathbb{E}\max_{j \in [\![1,k]\!]} \left | \left\langle Z, m(Q_{t_j,h})-(Y+Z)\right\rangle \right | \\
                   & \leq \sigma^2 + 2\sigma (\max_{j \in [\![1,k]\!]}\|m(Q_{t_j,h})\|+ \sqrt{2}(K+\sigma)) \\
                   & \leq \frac{C\sigma K}{\sqrt{h}},
\end{align*}
using \eqref{eq:bounds_sg_localmean} and $\sigma \leq K$. The converse bound on $P\dd^2_{Q,h,k} - Q\dd^2_{Q,h,k}$ may be proved the same way. Similarly, we may write
\begin{align*}
P\dd^2_{Q_n,h,k} - Q\dd^2_{Q_n,h,k} & \leq \sigma^2 + 2\sigma (\max_{j \in [\![1,k]\!]}\|m(Q_{n,t_j,h})\|+ \sqrt{2}(K+\sigma) \\
  & \leq \sigma^2 + 2\sigma(\max_{j \in [\![1,k]\!]}\|m(Q_{t_j,h})\| + \frac{C(K+\sigma)(p+1)\log(n)}{h\sqrt{n}}) +\sqrt{2}(K+\sigma)) \\
  & \leq \frac{C\sigma K(p+1)\log(n)}{h\sqrt{n}},
\end{align*} 
according to Lemma \ref{lm_ m bounded for subgaussian}. The bound on $Q\dd^2_{Q_n,h,k} - P\dd^2_{Q_n,h,k}$ derives from the same argument. Collecting all pieces, we have
\begin{align*}
\left |P (d^2_{Q_n,h,k} - d^2_{Q,h,k}) \right | & \leq  \left |Q (d^2_{Q_n,h,k} - d^2_{Q,h,k}) \right | + \frac{C\sigma K(p+1)\log(n)}{h\sqrt{n}} \\
               & \leq \frac{C\sigma K(p+1)\log(n)}{h\sqrt{n}} + \frac{CkK^2 ((p+1)\log(n))^{\frac{3}{2}}}{h \sqrt{n}},
\end{align*}
where we used $\sigma \leq K$.

\epv

\subsection{Proof of Lemma \ref{lem:Borne Q-Qn}}\label{sec:proof_lemma_Q_Qn}
\pv
With the notation $l_{t_i}(x)=-2\langle x,m(Q_{t_i,h})\rangle+\|m(Q_{t_i,h})\|^2+v(Q_{t_i,h})$, we get that:
\[
\sup_{t\in\R_d^{(k)}}(Q-Q_n)\gamma(t,u) = \sup_{t\in\R_d^{(k)}}\left(\left(Q-Q_n\right)\min_{i\in[\![1,k]\!]}l_{t_i}(u)\right).\]
First we note that since $Q$ is sub-Gaussian with variance $V^2$, we have, for every $c \in \overline{\R}_d$,
\begin{align}\label{eq:bounds_sg_localmean}
\|m(Q_{c,h})\|^2 + v(Q_{c,h}) = Q_{t,h}\|u\|^2 \leq \frac{2V^2}{h}.
\end{align}
Set $z=2V\sqrt{\log(n) + \lambda}$ and $\lambda= p\log(n)$. Then, with probability larger than $1-n^{-\frac{1}{p}}$, 
\begin{align}\label{eq:bound_sg_supradius}
\max_{i=1, \hdots, n}\|X_i\| \leq z.
\end{align}
We may then write
\begin{align*}
\sup_{t\in\R_d^{(k)}} \left | (Q-Q_n)\gamma(t,u) \right | & =   \sup_{t\in\R_d^{(k)}} \left | \frac{1}{n} \sum_{i=1}^{n} \gamma(t,X_i) - \mathbb{E}(\gamma(t,X)) \right | \\
                                 & \leq \sup_{t\in\R_d^{(k)}} \left | \frac{1}{n} \sum_{i=1}^{n} \gamma(t,X_i)\1_{\|X_i\|\leq z} - \mathbb{E}(\gamma(t,X) \1_{\|X\|\leq z)}) \right | \\
                                 & \quad + \sup_{t\in\R_d^{(k)}} \E( | \gamma(t,X)\1_{\|X\|>z} |) + \sup_{t\in\R_d^{(k)}} \left | \frac{1}{n} \sum_{i=1}^{n} \gamma(t,X_i) \1_{\|X_i\|>z} \right |.
\end{align*}
According to \eqref{eq:bound_sg_supradius}, the last part is $0$ with probability larger than $1-n^{-\frac{1}{p}}$. Moreover
\begin{align*}
\E(|\gamma(t,X)|\1_{\|X\|>z}) & \leq \E ( \1_{\|X\|>z}\sup_{j=1, \hdots, k}2 |\langle x,m(Q_{t_j,h}) |\rangle|+\|m(Q_{t_j,h})\|^2+v(Q_{t_j,h}) \\
  & \leq \frac{2V^2}{h} \mathbb{P}(\|X\| >z) + 2\sqrt{2}\frac{V}{\sqrt{h}} \E(\|X\|\1_{\|X\|>z}) \\
  & \leq 10 \frac{V^2}{h} e^{-\frac{z^2}{2V^2}} \\
  & \leq 10 \frac{V^2}{h} n^{-(p+1)}.
\end{align*}
It remains to bound 
\[
\sup_{t\in\R_d^{(k)}} \left | (Q-Q_n)\gamma(t,u)\1_{\|u\| \leq z} \right |.
\]
Since for every $t$ and $u$, $|\gamma(t,u)\1_{\|u\| \leq z}| \leq (z + \frac{V\sqrt{2}}{\sqrt{h}})^2:=R$, \cite[Theorem 6.2]{Massart16} entails
\[
\mathbb{P} \left ( \sup_{t\in\R_d^{(k)}} \left | (Q-Q_n)\gamma(t,u)\1_{\|u\| \leq z} \right | \geq \mathbb{E} \sup_{t\in\R_d^{(k)}} \left | (Q-Q_n)\gamma(t,u)\1_{\|u\| \leq z} \right | + R \sqrt{\frac{2\lambda}{n}} \right ) \leq e^{-\lambda} = n^{-p}.
\]
To bound $\mathbb{E} \sup_{t\in\R_d^{(k)}} \left | (Q-Q_n)\gamma(t,u)\1_{\|u\| \leq z} \right |$, we follow the same line as for Lemma \ref{lem_concentrationmoy}. A symmetrization argument yields
\begin{align*}\E \sup_{t\in\R_d^{(k)}} \left |(Q-Q_n)\min_{i\in[\![1,k]\!]}l_{t_i}(u)\1_{\|u\| \leq z} \right | & \leq \frac{2}{n}\mathbb{E}_{X_{1:n}} \mathbb{E}_{\sigma}\left[\sup_{t\in\R_d^{(k)}}\sum_{i=1}^n\sigma_i\min_{j\in[\![1,k]\!]}l_{t_j}(X_i)\1_{\|X_i\| \leq z}\right] \\
& \leq \frac{2R}{n}\mathbb{E}_{X_{1:n}} \mathbb{E}_{\sigma}\left[\sup_{t\in\R_d^{(k)}}\sum_{i=1}^n\sigma_i\min_{j\in[\![1,k]\!]}\frac{l_{t_j}(X_i)\1_{\|X_i\| \leq z}}{R}\right],
\end{align*}
where the $\sigma_i$'s are i.i.d. Rademacher variables, independent of the $X_i$'s, and $\mathbb{E}_Y$ denotes expectation with respect to the random variable $Y$. As in Section \ref{sec:proof_lemma_concentrationmoy}, denote, for any subset of functions $G$, $\mathcal{N}( G,\varepsilon)$ the $\varepsilon$-covering number of $G$ with respect to the metric $L_2(P_n)$. Denote by $\mathcal{F}_k$ the set of functions $ x \mapsto \min_{j\in[\![1,k]\!]}\frac{l_{t_j}(x)\1_{\|x\| \leq z}}{R}$, and by $\mathcal{F}$ the set of functions $x \mapsto \frac{l_{t}(x)\1_{\|x\| \leq z}}{R}$, $t \in \bar{\mathbb{R}}^d$.  Since the $x \mapsto\frac{l_{t_j}(x)\1_{\|x\| \leq z}}{R}$ are bounded by $1$, we may write, for any $\varepsilon>0$, 
\begin{align*}
\mathcal{N}(\mathcal{F}_k,\varepsilon) \leq \mathcal{N}(\mathcal{F},\varepsilon)^k,
\end{align*}
as well as
\begin{align*}
\mathcal{N}(\mathcal{F},\varepsilon) & \leq \mathcal{N}\left ( \left \{ x \mapsto \frac{-2 \left\langle x, m(Q_{t,h}) \right\rangle \1_{\|x\| \leq z}}{R} \right \}, \varepsilon/2\sqrt{2} \right ) \\ 
         & \qquad \qquad  \times \mathcal{N}\left ( \left \{ x \mapsto \frac{(\|m(Q_{t,h})\|^2+v(Q_{t,h})) \1_{\|x\| \leq z}}{R} \right \}, \varepsilon/2\sqrt{2} \right ) \\
         & \leq \mathcal{N}(\mathcal{G}_1,\varepsilon/2\sqrt{2}) \times \mathcal{N}(\mathcal{G}_2,\varepsilon/2\sqrt{2}).
\end{align*}
Using \cite[Theorem 1]{Mendelson03} yields
\begin{align*}
\mathcal{N}(\mathcal{G}_1, \varepsilon)) & \leq \left ( \frac{2}{\varepsilon} \right )^{2(d+1)C} \\
\mathcal{N}(\mathcal{G}_2,\varepsilon) & \leq \left ( \frac{2}{\varepsilon} \right )^{2C}.
\end{align*}
We then deduce
\begin{align*}
\mathcal{N} \left ( \mathcal{F}_k, \varepsilon \right ) \leq \left ( \frac{4 \sqrt{2}}{\varepsilon} \right )^{2k(d+2)C}.
\end{align*}
Using \cite[Corollary 13.2]{Massart16} leads to 
\begin{align*}
\mathbb{E}_{\sigma}\left[\sup_{t\in\R_d^{(k)}}\sum_{i=1}^n\sigma_i\min_{j\in[\![1,k]\!]}\frac{l_{t_j}(X_i)\1_{\|X_i\| \leq z}}{R}\right] \leq C \sqrt{k(d+2)n}.
\end{align*}
Combining these bounds gives the result of Lemma \ref{lem:Borne Q-Qn}.
\epv
\subsection{Proof of Lemma \ref{Borne gamma-hat gamma}}\label{sec:proof_borne_gammahat_gamma}
\pv
For $t\in\R_d^{(k)}$, we get that:
\begin{align*}
|\gamma(t,x)-\hat\gamma(t,x)|& \leq\max_{j\in[\![1,k]\!]} |-2\langle x,m(Q_{t_j,h})-m(Q_{n\,t_j,h})\rangle+\left(M(Q_{t_j,h})-M(Q_{n\,t_j,h})\right)| \\
& \leq 2\|x\| \max_{j\in[\![1,k]\!]} \|m(Q_{t_j,h})-m(Q_{n\,t_j,h}\| +  \max_{j\in[\![1,k]\!]}\|M(Q_{t_j,h})-M(Q_{n\,t_j,h}\|.
\end{align*}

Let $t \in \mathbb{R}^d$, and denote by $r=\delta_{Q,h}(t)$, $r_n=\delta_{Q_n,h}(t)$, and $z=2V\sqrt{(p+1)\log(n)}$. We may write
\begin{align*}
\|m(Q_{t,h})-m(Q_{n\,t,h})\|&\leq\frac{1}{h}\left(\|Qu\1_{\B(t,r)}(u)-Qu\1_{\B(t,r_n)}(u)\|+\|Qu\1_{\B(t,r_n)}(u)-Q_nu\1_{\B(t,r_n)}(u)\|\right)\\
&\leq\frac{1}{h}\left(\|(Q-Q_n)u\1_{\B(t,r_n)}(u)\| + Q\|u\||\1_{\B(t,r_n)} - \1_{\B(t,r)}|(u) \right) \\
&\leq\frac{1}{h}\left(\|(Q-Q_n)u\1_{\B(t,r_n)}(u)\| + z Q|\1_{\B(t,r_n)} - \1_{\B(t,r)}|(u) + Q\|u\|\1_{\|u\|>z} \right).
\end{align*}
Moreover, $Q\left|\1_{\B(t,r)}-\1_{\B(t,r_n)}\right|(u)=|h-Q(\B(t,r_n))|=|Q_n(\B(t,r_n))-Q(\B(t,r_n))|$. On the event described in Lemma \ref{lem_concentrationmoy}, we have that
\begin{align*}
\|(Q-Q_n)u\1_{\B(t,r_n)}(u)\| & \leq  C V \sqrt{d} \frac{(p+1)\log(n)}{\sqrt{n}}, \\
|Q_n(\B(t,r_n))-Q(\B(t,r_n))| & \leq C \sqrt{d} \frac{\sqrt{(p+1)\log(n)}}{\sqrt{n}}, \\
Q\|u\|\1_{\|u\|>z} & \leq 2Vn^{-(p+1)}.
\end{align*}
Thus,
\[\sup_{t\in\R^d}\|m(Q_{t,h})-m(Q_{n\,t,h})\|\leq \frac{CV (p+1) \log(n)}{h\sqrt{n}}.
\]
As well,
\begin{align*}
\sup_{t\in\R^d}|M(Q_{t,h})-M(Q_{n\,t,h})|& \leq\frac{1}{h}\left [ \left |(Q-Q_n) \|u\|^2\1_{\B(t,r_n)} \right | + Q \|u\|^2 | \1_{\B(t,r_n)} - \1_{\B(t,r}| \right ] \\
& \leq \frac{1}{h}\left [ \left |(Q-Q_n) \|u\|^2\1_{\B(t,r_n)} \right | + Q \|u\|^2\1_{\|u\|>z} + z^2 |(Q-Q_n)\1_{B(t,r_n)}| \right ]. 
\end{align*}
Using Lemma \ref{lem_concentrationmoy} again, we get
\begin{align*}
\left |(Q-Q_n) \|u\|^2\1_{\B(t,r_n)} \right | & \leq C V^2 \sqrt{d} \frac{(p+1)\log(n)^{\frac{3}{2}}}{\sqrt{n}} \\
|(Q-Q_n)\1_{B(t,r_n)}| & \leq C\sqrt{d} \frac{\sqrt{(p+1)\log(n)}}{\sqrt{n}} \\
Q \|u\|^2\1_{\|u\|>z} & \leq 2V^2n^{-(p+1)}.
\end{align*}
Collecting all pieces leads to
\begin{align}\label{eq:gamma_hatgamma}
\left |\gamma(t,x) - \hat{\gamma}(t,x)\right | \leq C \|x\|\frac{V (p+1) \log(n)}{h\sqrt{n}} + C V^2 \frac{(p+1)^{\frac{3}{2}}\log(n)^{\frac{3}{2}}}{h \sqrt{n}}.
\end{align}
At last, since 
\[
           \mathbb{P} \left \{ \max_{i} \|X_i \| \geq z \right \} \leq n e^{-\frac{t^2}{2V^2}} \leq n^{-2p+1},
           \]
we deduce that
\[
Q_n   \left |\gamma(t,x) - \hat{\gamma}(t,x)\right | \leq         C V^2 \frac{(p+1)^{\frac{3}{2}}\log(n)^{\frac{3}{2}}}{h \sqrt{n}}.
\]
\epv

\end{appendix}
\end{document}